\documentclass[a4paper, reqno]{amsart}
\usepackage[all]{xy}
\usepackage{amssymb,amsmath,amsthm}
\numberwithin{equation}{section}

\newcommand{\Q}{\mathbb{Q}}

\newcommand{\R}{\mathbb{R}}
\newcommand{\C}{\mathbb{C}}
\newcommand{\N}{\mathbb{N}}
\newcommand{\Z}{\mathbb{Z}}

\newcommand{\cQ}{\overline{\mathbb{Q}}}

\newtheorem{thm}{Theorem}
\newtheorem{lem}{Lemma}

\renewcommand{\mod}[1]{\hspace{-2.9mm}\pmod{#1}}
\newcommand{\x}{{\bf x}}

\newcommand{\bfP}{\mathbb{P}}

\newcommand{\ma}{\mathbf}
\newcommand{\ben}{\begin{enumerate}}
\newcommand{\een}{\end{enumerate}}
\newcommand{\eit}{\begin{itemize}}

\newcommand{\ve}{\varepsilon}

\newcommand{\mcal}{\mathcal}
\newcommand{\lab}{\label}
\newcommand{\al}{\alpha}

\newcommand{\D}{\Delta}
\newcommand{\be}{\beta}
\newcommand{\la}{\lambda}

\newcommand{\colt}[2]{\genfrac{}{}{0pt}{1}{#1}{#2}}

\newcommand{\tX}{{\widetilde X}}
\DeclareMathOperator{\hcf}{gcd}
\DeclareMathOperator{\meas}{meas}
\DeclareMathOperator{\pic}{Pic}

\newcommand{\vr}{\varrho}
\newcommand{\vt}{\vartheta}
\renewcommand{\d}{\mathrm{d}}
\renewcommand{\leq}{\leqslant}
\renewcommand{\geq}{\geqslant}

\newcommand{\m}{\mathfrak{m}}

\theoremstyle{definition}
\newtheorem*{ack}{Acknowledgements}

\begin{document}

\title[On Manin's conjecture for singular del Pezzo surfaces]{On Manin's conjecture for singular del Pezzo surfaces of degree
     four, II}

\author{R. de la Bret\`eche}
\author{T.D. Browning}

\address{
Institut de Math\'ematiques de Jussieu,
Universit\'e Paris 7 Denis Diderot,
Case Postale 7012,
2, Place Jussieu, 
F-75251 Paris cedex 05}
\email{breteche@math.jussieu.fr}

\address{School of Mathematics,  
University of Bristol, Bristol BS8 1TW}
\email{t.d.browning@bristol.ac.uk}


\begin{abstract}
This paper establishes the Manin conjecture for a 
certain non-split singular del Pezzo surface of degree four $X \subset \bfP^4$.
In fact, if $U \subset X$ is the open subset formed by deleting the
lines from $X$, and $H$ is the usual
projective height function on $\bfP^4(\Q)$, then the height zeta function
$
\sum_{x \in U(\Q)}{H(x)^{-s}}
$
is analytically continued to the half-plane $\Re e (s)>17/20$.
\end{abstract}

\maketitle

\section{Introduction}

Let $X \subset \bfP^4$ be a singular del Pezzo surface of degree four
such that $X(\Q)$ is Zariski dense in $X$,
and let $U \subset X$ denote the open subset formed by deleting the
lines from $X$.  The purpose of this paper is to extend our previous
investigation \cite{1} into the asymptotic distribution of rational points on $U$.
For any $x=[x_0,\ldots,x_4] \in \bfP^4(\Q)$ such that
$x_0,\ldots,x_4 \in \Z$ and $\hcf(x_0,\ldots,x_4)=1$, let
$H(x)=\max_{0\leq i \leq 4}|x_i|$
denote the usual anticanonical height function.
Then the behaviour of the associated counting function
$$
N_{U,H}(B)=\#\{x \in U(\Q): H(x) \leq B\},
$$
as $B \rightarrow \infty$, is predicted by the Manin conjecture
\cite{f-m-t}. 
Let $\tX$ denote the minimal 
desingularisation of $X$ and let $\rho$ denote the rank of the Picard group
$\pic \tX$ of $\tX$.  There is a strong version of this
conjecture that predicts the existence of a
constant $c_{X,H}>0$, and a monic polynomial $P\in \R[t]$ of
degree $\rho-1$, such that
\begin{equation}\lab{manin}
N_{U,H}(B)=c_{X,H} B P(\log B) +O(B^{1-\delta}),
\end{equation}
for some $\delta>0$.  
The constant $c_{X,H}$ has received a conjectural interpretation
at the hands of Peyre \cite{p}, and will be discussed in \S
\ref{conform} below.  As yet there appears to be no conjectural
understanding of the lower order coefficients in this asymptotic
formula. The true nature of the error term has been investigated
by Swinnerton-Dyer \cite{swd}, in the setting of diagonal cubic surfaces.

A classification of singular quartic del Pezzo surfaces  can
be found in the work of Hodge and Pedoe \cite[Book IV, \S
XIII.11]{h-p}. This shows that up to isomorphism over $\overline{\Q}$, there are
$15$ possible singularity types that can occur.
Coray and Tsfasman  \cite[Proposition 6.1]{c-t} have calculated the
extended Dynkin diagrams for each type.
Given this finite list of surfaces, it is natural to try and develop
an arsenal of tools and techniques that permit us to verify the
conjectured asymptotic formula \eqref{manin} for each surface on the list.
One approach to establishing the Manin conjecture involves studying
the height zeta function
$$
Z_{U,H}(s)=\sum_{x \in U(\Q)}\frac{1}{H(x)^s}.
$$
This is defined when $\Re e (s)$ is sufficiently large. 
Once one has proved suitably
strong statements about the analytic properties of $Z_{U,H}(s)$,  
one automatically obtains information about
the asymptotic behaviour of $N_{U,H}(B)$ via standard Tauberian arguments.
This approach was present in our previous work \cite{1}, where an 
extensive study was made of the quartic del Pezzo surface
\begin{equation}\lab{surface1}
x_0x_1-x_2^2=x_0x_4-x_1x_2+x_3^2 = 0.
\end{equation}
This surface is split over $\Q$ and has a unique singular point,
which is of type ${\mathbf D}_5$. In particular, the
Picard group of the minimal desingularisation of \eqref{surface1} has
maximal rank $6$.  In addition to providing an analytic continuation
of the corresponding height zeta function to the half-plane $\Re e
(s)>9/10$, an estimate of the shape \eqref{manin} was obtained for 
any $\delta \in (0,1/12)$.

The primary goal of this paper is to determine
whether the techniques that were developed in 
\cite{1} can be brought to bear upon a surface that is not split over
the ground field $\Q$. 
Let $X \subset \bfP^4$ be the surface 
\begin{equation}\lab{surface2}
x_0x_1-x_2^2=x_0^2-x_1x_4+x_3^2= 0.
\end{equation}
Then $X$ has a unique singular point $\xi=[0,0,0,0,1]$, which is of
type ${\mathbf D}_4$.
In fact $X$ has singularity 
type $\mathbf{C}_3$ over $\Q$, in the sense of Lipman \cite[\S
24]{lipman}, which becomes a $\mathbf{D}_4$ singularity over $\overline{\Q}$.
It is easy to see that the only two lines that are contained in $X$  are
$$
\ell_1: x_1=x_2=x_0- i x_3=0, \qquad \ell_2: x_1=x_2=x_0+ i x_3=0.
$$
Clearly both $\ell_1$ and $\ell_2$ pass through $\xi,$ which is actually
the only rational point lying on either line.
In particular
$N_{U,H}(B)=N_{X,H}(B)+O(1).$
In our previous work \cite{1} the universal torsor was a fundamental
ingredient in the resolution of the 
Manin conjecture for \eqref{surface1}, much 
in keeping with the general philosophy.
One of the most novel features of our present investigation is that
we will be able to establish the Manin conjecture for
\eqref{surface2} using a certain sub-torsor of the universal torsor.

Our first result concerns 
the analytic properties of the associated height zeta function $Z_{U,H}(s)$.
For any positive integer $n$, let
$$
\chi(n)=
\left\{
\begin{array}{ll}
+1  & \hbox{if } n\equiv 1  \mod{4},\\
-1  & \hbox{if } n\equiv 3  \mod{4},\\
0  & \hbox{otherwise},
\end{array}
\right.
$$
be the real non-principal character modulo $4$.
Then for $\Re e(s)>0$ we introduce the functions
\begin{align}
E_1(s+1)&=\zeta(2s+1)^2\zeta(3s+1)\zeta(4s+1)
L(2s+1,\chi)L(3s+1,\chi),\lab{e1}\\
E_2(s+1)&=\frac{\zeta(9s+3)L(9s+3,\chi)}{\zeta(5s+2)^2
\zeta(6s+2)^2L(5s+2,\chi)L(6s+2,\chi)^2}.\lab{e2}
\end{align}
It is easily seen that $E_1(s)$ has a meromorphic
analytic continuation to the entire complex plane with a single pole at $s=1$.
Similarly it is clear that $E_2(s)$ is holomorphic and bounded
on the half-plane $\{ s\in\C: \Re e (s)\geq 5/6+\ve\}$, for any $\ve>0$.
For any $\al\in \R$ let 
$$
\mcal{H}_\al=
\{ s\in\C: \Re e (s)\geq \al+\ve\}.
$$
We are now ready to state our main result.

\begin{thm}\lab{main'}
Let $\ve>0$.  Then there exists a constant $\beta
\in \R$, and functions $G_1(s), G_2(s)$
that are holomorphic on the half-planes $\mcal{H}_{3/4}$ and
$\mcal{H}_{17/20}$, respectively, 
such that for $\Re e(s)>1$ we have
$$
Z_{U,H}(s) =E_1(s)E_2(s)G_1(s) +\frac{{12/\pi^2}+4\be}{s-1}+G_2(s).
$$
In particular $(s-1)^4Z_{U,H}(s)$ has a holomorphic analytic
continuation to the half-plane $\mcal{H}_{17/20}$.
The function $G_1(s)$ is bounded on the half-plane $\mcal{H}_{3/4}$
and satisfies $G_1(1)\neq 0$, 
and the function $G_2(s)$ satisfies  
$$ 
G_2(s)\ll_\ve (1+|\Im m (s)|)^{20\max\{1-\Re e (s), 0\}/3+\varepsilon} 
$$
on the domain $\mcal{H}_{17/20}$.
\end{thm}

The main step in the proof of Theorem
\ref{main'} consists of establishing
a preliminary estimate for $N_{U,H}(B)$.  This will be the object of
\S\S \ref{pre-man}--\ref{final}.
In \S \ref{height} this estimate
will then be used to deduce the analytic properties of $Z_{U,H}(s)$ 
presented above.
Explicit expressions for $\beta$ and $G_1$ can be found in
\eqref{beta} and \eqref{calculG1}, respectively.
It is interesting to compare Theorem \ref{main'} with the
corresponding result in our previous work \cite[Theorem 1]{1}.
It is no surprise that the structure of the two height zeta functions
is very similar.  Thus in both expressions we have a first
term $E_1(s)E_2(s)G_1(s)$ that corresponds to the main term in our 
preliminary estimate
for the counting function, a term $\frac{12}{\pi^2}(s-1)^{-1}$ 
that corresponds to an isolated conic contained in the surface, and
a further ``$\beta$-term'' involving a constant $\beta$ that arises  
through the error in approximating certain arithmetic 
quantities by real-valued
continuous functions. This $\beta$-term is one of the most mysterious
aspects of our work, and it is interesting to highlight the difference
in nature between the constant that appears in Theorem \ref{main'} and
the corresponding constant obtained in \cite[Theorem 1]{1}.  
Thus whereas the $\beta$-term in the latter work  
relies upon results concerning the equidistribution of squares in a
fixed residue class, the $\beta$-term in Theorem \ref{main'} merely 
arises through a routine application of integration by
parts.  

We have already mentioned that the main step in the proof of Theorem
\ref{main'} involves producing a preliminary estimate for
$N_{U,H}(B)$. In \S \ref{deduc} we will show how Perron's formula 
can be combined with Theorem~\ref{main'} to 
extract the following asymptotic formula for $N_{U,H}(B)$.

\begin{thm}\lab{main}
Let $\delta \in (0,3/32)$.
Then there exists a polynomial $P$ of
degree $3$ such that  for any $B \geq 1$ we have
$$
N_{U,H}(B)=B P(\log B) +O(B^{1-\delta}).
$$
Moreover, the leading coefficient of $P$ is equal to
$$
\frac{\pi^2}{ 576}\int_0^1 \frac{u^{1/4}\d u}{ \sqrt{1-u}}
\prod_{p }\Big(1-\frac{1}{p}\Big)^{4}\Big(1-\frac{\chi(p)}{p }\Big)^2
\Big(1+\frac{3+2\chi(p)+\chi^2(p)}{p}+ \frac{\chi^2(p)}{p^2}\Big).
$$
\end{thm}

We will verify in \S \ref{conform}
that Theorem \ref{main} is in accordance with Manin's conjecture.
A crucial step in the proof of Theorems \ref{main'} and \ref{main} is 
a bijection that we establish between the rational points on
$U$ and the points $(v_1,v_2,y_0,\ldots,y_4) \in \Z^7$ such that
\begin{equation}\lab{torsor}
y_0^4y_2^2   -v_2y_1^2y_4+ y_3^2=0.
\end{equation}
Note that $v_1$ does not appear explicitly in the equation.
This step is achieved in \S \ref{pre-man} via an elementary analysis of
the equations defining $X$.  As we have already indicated 
it is interesting to note that this
equation is not an affine embedding of the universal torsor over the minimal desingularisation
$\tX$ of $X$.  Instead it turns out that \eqref{torsor} corresponds to a
certain sub-torsor of the universal torsor, which reflects the fact
that $X$ does not split over the ground field.
Theorem \ref{main} seems to signify the first time that the full Manin
conjecture has been established without recourse to the universal torsor.

Over the last decade or so the Manin conjecture has been established
for a variety of special cases, and it is important to place 
our investigation in the context of other work.
We will say nothing about the situation for 
non-singular del Pezzo surfaces, or singular del Pezzo 
surfaces of degree not equal to four.  A discussion of these surfaces can be found
in the second author's survey \cite{gauss}.
It turns out that the Manin conjecture has already been established
for several singular quartic del
Pezzo surfaces by virtue of the fact that the surface is toric, for
which there is the general work of Batyrev and Tschinkel \cite{b-t'},
or the surface is an equivariant compactification of 
$\mathbb{G}_a^2$, for which there is the work of Chambert-Loir and Tschinkel \cite{ct}.
The surface \eqref{surface1} studied in \cite{1}
falls into this latter category.  That \eqref{surface2} is not an equivariant
compactification of $\mathbb{G}_a^2$ can be seen by mimicking the
argument used by Hassett and Tschinkel \cite[Remark 4.3]{ht} in their
analysis of a certain cubic surface.
The authors have recently learnt of work due to Derenthal and
Tschinkel \cite{d-t}, in which the Manin conjecture is established for
the surface
$$
x_0x_3-x_1x_4=x_0x_1+x_1x_3+x_2^2= 0.
$$
This is the split del Pezzo surface of degree four, with singularity type ${\mathbf D}_4$.
Their asymptotic formula is weaker than ours, and does not lead to an
analytic continuation of the corresponding height zeta function.


\begin{ack}
The authors are grateful to Ulrich
Derenthal and Brendan Hassett for several useful
conversations relating to universal torsors for singular del Pezzo surfaces.
Special thanks are due to Roger Heath-Brown whose ideas led
us to the proof of Lemma \ref{av-order}.
The paper was finalised while the first author was
at the \'Ecole Normale Sup\'erieure, 
and the second author was at Oxford University supported by 
EPSRC grant number GR/R93155/01.
The hospitality and financial support of these institutions is
gratefully acknowledged. 
Finally, the authors would like to thank the anonymous
referee for his careful reading of the manuscript and numerous useful suggestions.
\end{ack}

\section{Conformity with the Manin conjecture}\lab{conform}

In this section we will review some of the geometry of the surface $X\subset
\bfP^4$, as defined by the pair of quadratic forms
\begin{equation}\lab{forms}
Q_1(\x)=x_0x_1-x_2^2, \quad Q_2(\x)=x_0^2-x_1x_4+x_3^2,
\end{equation}
where $\x=(x_0,x_1,x_2,x_3,x_4)$.
In particular we will show that Theorem \ref{main} agrees with the
Manin conjecture.

Let $\tX$ denote the minimal desingularisation
of $X$,  and let $\pi: \tX\rightarrow X$ denote the corresponding
blow-up map.    We let $L_i$ denote the strict transform of the line
$\ell_i$ for $i=1,2$, and let $E_1,\ldots,E_4$ denote the exceptional
curves of $\pi$.  Then the divisors $E_1,\ldots,E_4, L_1, L_2$
satisfy the Dynkin diagram
$$
\xymatrix{
        & & E_2 \ar@{-}[d] \\
L_1 \ar@{-}[r] & E_3 \ar@{-}[r] & E_1 \ar@{-}[r] & E_4
\ar@{-}[r] & L_2 }
$$
after a possible relabelling of indices. From this it is possible to
write down the $6\times 6$ intersection matrix 
\begin{center}
\begin{tabular}{c|rrrrrr}
       & $E_1$ & $E_2$ & $E_3$ & $E_4$& $L_1$ & $L_2$\\
\hline
$E_1$ & $-2$ & $1$ & $1$ & $1$ & $0$ & $0$  \\
$E_2$ & $1$ & $-2$ & $0$ & $0$ & $0$ & $0$  \\
$E_3$ & $1$ & $0$ & $-2$ & $0$ & $1$ & $0$  \\
$E_4$ & $1$ & $0$ & $0$ & $-2$ & $0$ & $1$  \\
$L_1$ & $0$ & $0$ & $1$ & $0$ & $-1$ & $0$  \\
$L_2$ & $0$ & $0$ & $0$ & $1$ & $0$ & $-1$
\end{tabular}
\end{center}
which implies that the geometric Picard
group $\pic_{\cQ} \tX=\pic(\tX \otimes_\Q \cQ)$ is generated by $E_1,E_2,E_3,E_4,L_1,L_2$.
Moreover the adjunction formula implies that
\begin{equation}\lab{21-K}
-K_{\tX}=4E_1+2E_2+3(E_3+E_4)+2(L_1+L_2),
\end{equation}
where $-K_{\tX}$ denotes the anticanonical divisor of $\tX$.
Now if $\Gamma$ denotes the Galois group of
$\Q(i)/\Q$, then it is clear that  $\{L_1,L_2\}^\sigma=\{L_1,L_2\}$
for any $\sigma \in \Gamma$.
Furthermore, it emerges during the calculation of $\tX$ that
$$
E_1^\sigma=E_1, \quad
E_2^\sigma=E_2, \quad \{E_3,E_4\}^\sigma=\{E_3,E_4\},
$$
for any $\sigma \in \Gamma$.  The Picard group $\pic \tX$ of
$\tX$ is therefore the  free abelian group generated by 
$E_1,E_2,E_3+E_4, L_1+L_2$. 
In particular $\rho=4$ in \eqref{manin}, which agrees
with Theorem \ref{main}.

It remains to discuss the conjectured value of the constant $c_{X,H}$ in
\eqref{manin}.  For this we will follow the presentation 
adopted in our previous investigation \cite[\S 2]{1}, and so we will permit
ourselves to be brief.  In the notation found there, we see that 
the conjectured value of the constant in \eqref{manin} is 
\begin{equation}\lab{constant}
c_{X,H}=\alpha(\tX)\beta(\tX)\tau_H(\tX).
\end{equation}
Now it follows from \cite[Theorem 7.2]{c-t} that
$H^1(\Q,\pic_{\cQ}{\tX})=0$, whence $\beta(\tX)=1$. 
Turning to the value of $\alpha(\tX)$, we have already seen how
$-K_{\tX}$ can be written in terms of the basis for $\pic\tX$.
Moreover it is easy to check that the cone of
effective divisors  $\Lambda_{\mathrm{eff}}(\tX)\subset (\pic{\tX}) \otimes_\Z
\R$ is also generated by the basis elements of $\pic
\tX$.  This allows us to conclude that
\begin{equation}\lab{alpha}
\begin{split}
\alpha(\tX)&= \meas\big\{(t_1,t_2,t_3,t_4)\in \R_{\geq 0}^4:
4t_1+2t_2+3t_3+2t_4=1 \big\}\\
&= \frac{1}{2 }
\meas\big\{(t_1,t_2,t_3)\in \R_{\geq 0}^3:
4t_1+2t_2+3t_3\leq 1 \big\}\\
&=\frac{1}{288}.
\end{split}
\end{equation}
The calculation of $\tau_H(\tX)$ is a little more involved, and will be
carried out in the following result.

\begin{lem}\lab{localdensities}
We have $\tau_H(\tX)= \pi^2\omega_\infty\tau/16$, where
\begin{equation}\lab{definfty}
\omega_\infty=8\int_0^1 \frac{u^{1/4}\d
u}{ \sqrt{1-u}}
\end{equation}
and
\begin{equation}\lab{deftauf}
\tau=
\prod_{p }\Big(1-\frac{1}{p}\Big)^{4
}\Big(1-\frac{\chi(p)}{p }\Big)^2
\Big(1+\frac{3+2\chi(p)+\chi^2(p)}{p}+
\frac{\chi^2(p)}{p^2}\Big).
\end{equation}
\end{lem}

\begin{proof}
Write $L_p(s,\pic_{\cQ}\tX)$ for the local factors of
$L(s,\pic_{\cQ}\tX)$.  Furthermore, let
$\omega_\infty$ denote the archimedean density of points on $X$, and let
$\omega_p$ denote the usual $p$-adic density of points on $X$, for
any prime $p$.  Then the Tamagawa measure is given by
$$
\tau_H(\tX)= \lim_{s\rightarrow 1}\big((s-1)^\rho
L(s,\pic_{\cQ} \tX)\big)
\omega_\infty\prod_p
\frac{\omega_p}{L_p(1,\pic_{\cQ} \tX)},
$$
where $\rho=4$ is the rank of $\pic\tX$.
Our first step is to note that
$$
L(s,\pic_{\cQ} \tX)=\zeta(s)^2\zeta_{\Q(i)}(s)^2.
$$
Since $\zeta_{\Q(i)}(s)=\zeta(s)L(s,\chi)$, it then follows that
$$
\lim_{s\rightarrow 1}\big((s-1)^\rho
L(s,\pic_{\cQ} \tX)\big)
=L(1,\chi)^2=\frac{\pi^2}{16}.
$$
Furthermore, we plainly  have
\begin{equation}\lab{Lp1}
L_p(1,\pic_{\cQ} \tX)^{-1}=
\Big(1-\frac{1}{p}\Big)^4\Big(1-\frac{\chi(p)}{p}\Big)^2,
\end{equation}
for any prime $p$.

We proceed by calculating the value of the
archimedean density $\omega_\infty$. Let $\|\x\|$ denote the norm
$\max_{0 \leq i \leq 4} |x_i|$ for any $\x \in \R^5$.
We will follow the method given by Peyre \cite{p} to compute
$\omega_\infty$.
It will be convenient to parametrise
the points via the choice of variables $x_0,x_1,x_4$,
for which we first observe that
$$
\det \left(
\begin{matrix}
\frac{\partial Q_1}{ \partial x_2}&\frac{\partial Q_2}{
\partial x_2} \cr\cr  \frac{\partial Q_1}{ \partial x_3}
&\frac{\partial Q_2}{ \partial x_3}
\end{matrix}
\right)=-4x_2x_3.
$$
Now in any real solution to the pair of equations $Q_1(\x)=Q_2(\x)=0$,
the components $x_0,x_1$ and $x_4$ must necessarily all share the same
sign.  Taking into account the fact that $\x$ and $-\x$ represent the
same point in $\bfP^4$, the archimedean density of points on $X$
is therefore equal to
$$
\omega_{\infty}  = 4
\int_{ \{\x\in\R_{>0}^5: \,Q_1(\x)=Q_2(\x)=0,\, \|\x\|\leq 1 \}}\omega_L( \x),
$$
where $\omega_L(\x)$ is the Leray form
$(4x_2x_3)^{-1}\d x_0\d x_1\d x_4$. It follows that
\begin{align*}
\omega_{\infty} &=\int\int\int_{\big\{\colt{x_0,x_1,x_4\in\R_{>0}:}{x_1,x_4\leq
1,\,x_0^2<x_1x_4}\big\}}\frac{\d x_0\d x_1\d x_4}{
\sqrt{x_0x_1(x_1x_4-x_0^2)}}\\
&=2\int \int_{\big\{\colt{x_0,x_1 \in\R_{>0}:}{x_1 \leq
1,\,x_0^2< x_1 }\big\}}\frac{\sqrt{ x_1 -x_0^2 }}{
x_0^{1/2}x_1^{3/2}}\d x_0\d x_1.
\end{align*}
The change of variables $u=x_0^{2}/ x_1 $ therefore yields
$$
\omega_{\infty}= \int_0^1\frac{\sqrt{1-u}}{u^{3/4}}\d
u\int_0^1\frac{\d x_1}{x_1^{3/4}}=
8 \int_0^1 \frac{u^{1/4}\d u}{ \sqrt{1-u}},$$
where we have carried out integration by parts to get the last
equality. This establishes \eqref{definfty}.

It remains to calculate the value of
$\omega_p=\lim_{r\to \infty}p^{-3r}N(p^r)$, for any prime $p$, where
$$
N(p^r)=\#\{\x \mod{p^r}: Q_1(\x)\equiv
Q_2(\x)\equiv 0 \mod{p^r}\}.
$$
Although this amounts to a routine calculation, the arguments needed
to handle our non-split surface are slightly more subtle, and we have decided to
present them in full.  To begin with we write
$x_0=p^{k_0}x_0'$ and $x_1=p^{k_1}x_1'$, with $p\nmid x_0'x_1'$.
Now we have $p^r \mid x_2^2$ if and only if $k_0+k_1\geq r$,
and there are at most $p^{r/2}$ square roots of zero modulo $p^r$.
When $ k_0+k_1< r$, it follows that $k_0+k_1$ must be even and we may write
$x_2=p^{(k_0+k_1)/2}x_2',$ with $p\nmid x_2'$ and
$$
x_0'x_1'-{x_2'}^2\equiv 0 \mod{p^{r-k_0-k_1}}.
$$
The number of possible choices for $x_0',x_1',x_2'$ is therefore
$$
h_p(r,k_0,k_1)=\left\{
\begin{array}{ll} \phi(p^{r-k_0})\phi(p^{r-(k_0+k_1)/2})p^{k_0} &
\mbox{if $k_0+k_1< r$},\\
O(p^{5r/2-k_0- k_1 })& \mbox{if $k_0+k_1\geq r$}.
\end{array}
\right.
$$

It remains to determine the number of solutions $x_3,x_4$ modulo $p^r$
such that
\begin{equation}\lab{equationx3x4}
p^{2k_0}{x_0'}^2-p^{ k_1}x_1'x_4+x_3^2\equiv 0 \mod{p^r}.
\end{equation}
In order to do so we distinguish between four basic cases:  either
$k_0+k_1<r$ and $2k_0=k_1$, or
$k_0+k_1<r$ and $2k_0<k_1$, or
$k_0+k_1<r$ and $2k_0>k_1$, or else $k_0+k_1\geq r$. For the first
three of these cases we must take care only to sum over values of
$k_0,k_1$ such that $k_0+k_1$ is even.
We will denote by $N_i(p^r)$ the contribution to $N(p^r)$ from the
$i$th case, for $1\leq i \leq 4$, so that
\begin{equation}\label{N_i}
N(p^r)=N_1(p^r)+N_2(p^r)+N_3(p^r)+N_4(p^r).
\end{equation}

We begin by calculating the value of $N_1(p^r)$.  For this we write
$x_3=p^{k_3}x_3'$, with $k_3=\min \{ r/2, k_0\}=k_0$.  The number of
possibilities for $x_3'$ is $p^{r-k_0}$, each one leading to precisely
$p^{2k_0}$ possible choices for $x_4$ via \eqref{equationx3x4}.  On
noting that $k_1$ is even, so that $k_0$ must be even, we may write
$k_0=2k_0'$, for $0 \leq k_0'<r/6$.  In this way we deduce that
$$
N_1(p^r)=\sum_{0\leq k_0'< r/6}p^{r+2k_0'}h_p(r,2k_0',4k_0')
= p^{3r}(1-1/p)+O(p^{3r-r/6-1}).
$$

Next we calculate $N_2(p^r)$.  As above we write
$x_3=p^{k_0}x_3'$, and consider the resulting congruence
$$
{x_0'}^2+{x_3'}^2\equiv p^{k_1-2k_0}x_1'x_4 \mod{p^{r-2k_0}}.
$$
Suppose first that $p$ is odd.  Then modulo
$p^{k_1-2k_0}$, there are  $1+\chi (p)$ choices for $x_3'$,
whence there are $(1+\chi (p))p^{r+k_0-k_1}$ possibilities for $x_3'$
overall.   But then $x_4$ satisfies a congruence
modulo $p^{r-k_1}$, and there are therefore $p^{ k_1}$ ways of
choosing $x_4$.  On summing these contributions over all the relevant
values of $k_0,k_1$, we obtain
\begin{align*}
N_2(p^r)&=\sum_{\colt{k_0+k_1<r, ~k_0,k_1\geq 0}
{2k_0<k_1, ~2\mid (k_0+k_1)}}
(1+\chi(p))p^{r+k_0}h_p(r,k_0,k_1)\\
&= 2(1+\chi (p))p^{3r-1}\big(1+o(1)\big),
\end{align*}
when $p$ is odd. When $p=2$ the only difference in this calculation is
that we must restrict ourselves to the case
$k_1-2k_0=1$, since there are no solutions to the congruence
${x_0'}^2+{x_3'}^2\equiv 0 (\bmod{ 2^{\nu}})$ for $\nu \geq 2$.
We therefore obtain
\begin{align*}
N_2(2^r)&=\sum_{\colt{k_0+k_1<r, ~k_0,k_1\geq 0}
{ k_1=2k_0+1, ~2\mid (k_0+k_1)}}
(1+\chi(2))2^{r+k_0}h_2(r,k_0,k_1)\\
&=   2^{3r-2}\big(1+o(1)\big).
\end{align*}

We proceed by calculating $N_3(p^r)$, for which we now write
$x_3=p^{k_3}x_3'$ with
$k_3=\min \{ r/2,\lceil  k_1/2\rceil\}=\lceil  k_1/2\rceil$.
Here, as throughout this work, $\lceil \al \rceil$ denotes the ceiling
function for any $\al \in \R$, whereas $[\al]$ will always denote
the integer part of $\al$.
With this in mind \eqref{equationx3x4} becomes
$$
p^{2k_0-k_1}{x_0'}^2 +p^{2 \lceil  k_1/2\rceil
-k_1}{x_3'}^2\equiv  x_1'x_4 \mod{p^{r-k_1}}.
$$
If $k_1$ is even, the number of choices for $x_3$ is
$p^{r-k_1/2}$ and this leads to $p^{k_1}$ choices for $x_4$.  Thus
there are $p^{r+k_1/2}$ possibilities for $x_3,x_4$ if $k_1$ is even.
If $k_1$ is odd, then $p$ divides $x_4$ and we find that there are
$p^{r+(k_1-1)/2}$
possibilities for $x_3$ and $x_4$.  Summing over the relevant values
of $k_0$ and $k_1$ we therefore obtain
\begin{align*}
N_3(p^r)
&=\sum_{\colt{k_0+k_1<r, ~k_0,k_1\geq 0}
{2k_0>k_1, ~2\mid (k_0+k_1)}}
p^{r+[k_1/2]}h_p(r,k_0,k_1)\\
&= p^{3r-2}(p-1)^2 \sum_{\colt{k_0+k_1<r, ~k_0,k_1\geq 0}
{2k_0>k_1, ~2\mid (k_0+k_1)}}
p^{[k_1/2]-k_1/2-k_0/2}.
\end{align*}
On splitting the summation into four different cases according to the
value of $k_1$ modulo $4$, a routine calculation
therefore yields
\begin{align*}
\frac{N_3(p^r)}{p^{3r-2}(p-1)^2}&=
\sum_{\colt{ k_0=2k_0',~k_1=4k_1' }{ k_0'\geq k_1'+1}}p^{-k_0' }+
\sum_{\colt{ k_0=2k_0',~k_1=4k_1'+2 }{ k_0'\geq k_1'+1}}
p^{-k_0'}\\&\quad
+\sum_{\colt{ k_0=2k_0'+1,~k_1=4k_1'+1 }{ k_0'\geq k_1' }}
     p^{-k_0'-1 }+\sum_{\colt{ k_0=2k_0'+1,~k_1=4k_1'+3 }{ k_0'\geq k_1'+1}}
     p^{-k_0' -1}\\
&= \frac{3p+1}{(p-1)^2}\big(1+o(1)\big),
\end{align*}
as $r \rightarrow \infty$.

Finally we calculate the value of $N_4(p^r)$.  In this case a
straightforward calculation shows that there are at most
$2p^{r-k_1+\min\{k_0,k_1/2\}}$ possibilities for $x_3$, each one
leading to at most $p^{k_1}$ choices for $x_4$.  We therefore deduce that
$$
N_4(p^r)\ll
\sum_{k_0+k_1\geq r} p^{r+\min\{k_0,k_1/2\}}h_p(r,k_0,k_1) \ll p^{3r-r/6}.
$$
We may now combine our estimates for $N_1(p^r),\ldots, N_4(p^r)$ into
\eqref{N_i}.  When $p$ is odd we therefore deduce that
$$
\omega_p=\lim_{r\to\infty}p^{-3r}N(p^r)=1+\frac{4+2\chi(p)}{p}+\frac{1}{p^2},
$$
whereas when $p=2$ we obtain
$$
\omega_2=\lim_{r\to\infty}2^{-3r}N(2^r)=\frac52.
$$
On combining this with \eqref{Lp1}, we therefore conclude the proof of
\eqref{deftauf}, and so complete the proof of Lemma \ref{localdensities}.
\end{proof}

We end this section by combining \eqref{alpha} and Lemma
\ref{localdensities} in \eqref{constant}, in order to deduce that the
conjectured value of the constant in \eqref{manin} agrees with 
the value of the leading coefficient in Theorem \ref{main}.

\section{Preliminary manoeuvres}\lab{pre-man}

In this section we will establish an alternative expression for
$N_{U,H}(B)$, for which we will follow the presentation of \cite[\S 4]{1}.
Let us begin by recalling the notation used there.
For any $n \geq 2$ we will let $Z^{n+1}$ denote the set of
primitive vectors in $\Z^{n+1}$, and similarly, we let $N^{n+1}$ denote the set  of primitive vectors in $\N^{n+1}$.
Moreover, let $\Z_*^{n+1}$ (resp. $Z_*^{n+1}$) denote the
set of vectors $\ma{v} \in \Z^{n+1}$ (resp. $\ma{v} \in Z^{n+1}$)
such that $v_0\cdots v_n \neq 0$.
Finally, we will henceforth follow common convention and allow the
small parameter $\ve>0$ to take different values at different
points of the argument.

If $x=[\x]\in X\cap \bfP^4(\Q)$ is represented by
the vector $\x \in Z^5$, then it is easy to see that
$$
H(x)=\max\{|x_1|,|x_4|\}.
$$
Recall the definition \eqref{forms} of the quadratic forms $Q_1$ and $Q_2$.
Then it follows that
$$
N_{U,H}(B)= \frac{1}{2}\#\Big\{\x\in Z^5:
\max\{|x_1|,|x_4|\}\leq B, ~Q_1(\x)=Q_2(\x)=0\Big\} +O(1),
$$
since $\x$ and $-\x$ represent the same point in $\bfP^4$.
Let
\begin{equation}\label{india}
N(Q_1,Q_2;B)=\#\{\x\in N^5: \max\{x_1,x_4\}\leq B, ~Q_1(\x)=Q_2(\x)=0\}.
\end{equation}
Then we proceed to establish the following basic result.

\begin{lem}\lab{Reduc1}
Let $B\geq 1$.  Then we have
$$
N_{U,H}(B)=4 N(Q_1,Q_2;B) + \frac{12}{\pi^2}B + O(B^{2/3}).
$$
\end{lem}

\begin{proof}
Let us consider the contribution to $N_{U,H}(B)$ from vectors
$\x \in Z^5$ which contain zero components.  We claim that
\begin{equation}\lab{n(b)}
N_{U,H}(B)= \frac{1}{2}\#S(Q_1,Q_2;B) + \frac{12}{\pi^2}B+O(B^{2/3}),
\end{equation}
where
$$
S(Q_1,Q_2;B)=\{\x\in Z_*^5: \max\{|x_1|,|x_4|\}\leq B,
~Q_1(\x)=Q_2(\x)=0\}.
$$
Suppose that $\x\in Z^5$ is a vector such that
$$
x_0=0, \quad |x_1|,|x_4|\leq B,
$$
and $Q_1(\x)=Q_2(\x)=0$.  Then it immediately follows from the
first of these equations that $x_2=0$, and the second equation
implies that $x_3^2=x_1x_4$.
Now either $x_3=0$, in which case there are at most $4$ 
possibilities for $\x$, or else we have $\x=\pm(0,a^2,0,\pm ab,b^2)$
for coprime $a, b \in \N$. 
Hence the overall contribution from this case is $12B/\pi^2+O(B^{1/2})$. 
  Suppose now that  $\x\in Z^5$ is a vector such
that
$$
x_1=0, \quad |x_4|\leq B,
$$
and $Q_1(\x)=Q_2(\x)=0$.  Then a straightforward analysis of these
equations reveals that in fact $\x=\pm (0,0,0,0,1)$.  
Next we note that any vector
with $x_2x_4=0$ must have $x_0x_1=0$.
But such vectors have already been handled above.
Finally, if $\x\in Z^5$ satisfies
$$
x_3=0, \quad |x_1|,|x_4|\leq B,
$$
and $Q_1(\x)=Q_2(\x)=0$, then we must have
$x_1x_4=x_0^2$ and  $x_2^2=x_0x_1$.
Since we are only interested in an upper bound it clearly suffices
to count non-zero
integers $x_0,x_2, x_4$, with modulus at most $B$,  such that
$\hcf(x_0,x_2,x_4)=1$ and
$x_2^2x_4=x_0^3$.  But then it follows that
$(x_0,x_2,x_4)=\pm (a^2b, a^3, b^3)$ for coprime
integers $a,b$, whence the overall contribution is $O(B^{2/3})$.
This completes the proof
of \eqref{n(b)}.

We now need to relate the cardinality of the set $S(Q_1,Q_2;B)$ to
the quantity $N(Q_1,Q_2;B)$. Write $S=S(Q_1,Q_2;B)$ for convenience, and define the four subsets
$$
S_{\pm,\pm}=\{\x \in
S: \pm x_2>0, \pm x_3>0\}.
$$
Then we clearly have a
disjoint union
$S=S_{-,-} \cup S_{-,+} \cup S_{+,-} \cup S_{+,+}$,
in which each of the four sets has equal size.  Hence it follows that
$\#S=4\#S_{+,+}.$
Now for any $\x \in S_{+,+}$, we must have that $x_0$ and $x_1$ are
both positive or both negative, since their product is a square.  Similarly
$x_1$ and $x_4$ both have the same sign, since their product is the
sum of two squares.  Hence either $x_0,x_1,x_4$ are all positive, or
they are all negative.  This therefore establishes the equality
$$
\#S_{+,+}=2N(Q_1,Q_2;B).
$$
Upon recalling that $\#S=4\#S_{+,+}$, and then inserting this into
\eqref{n(b)}, we thereby complete the proof of Lemma \ref{Reduc1}.
\end{proof}

We proceed to equate $N(Q_1,Q_2;B)$ to the number of integral points
on a certain affine variety related to $X$,
subject to certain constraints.
Let $\x \in N^5$ be any vector counted by $N(Q_1,Q_2;B)$.  We begin
by considering solutions to the equation
$Q_1(\x)=0.$
But it is easy to see that there is a bijection
between the set of positive integers $x_0,x_1,x_2$ such that
$x_0x_1=x_2^2$, and the set of $x_0,x_1,x_2$ such that
$$
x_0=z_0^2z_2, \quad x_1=z_1^2z_2, \quad x_2=z_0z_1z_2,
$$
for $z_0, z_1,z_2 \in \N$, with
\begin{equation}\lab{polo-1}
\hcf(z_0,z_1)=1.
\end{equation}
We now substitute these values into the equation
$Q_2(\x)=0,$
in order to obtain
\begin{equation}\lab{golf-1}
z_0^4z_2^2-x_4z_1^2z_2+x_3^2=0.
\end{equation}
It is apparent that $z_2$ divides $x_3^2$.  Hence we write
$$
z_2=v_2{y_2'}^2,
$$
for $v_2 ,y_2' \in \N$ such that $v_2$ is square-free.  It follows
that $v_2 y_2' $ divides $x_3$, and so there exists $y_3' \in \N$ such that
$$
x_3=v_2 y_2' y_3'.
$$
Under these substitutions \eqref{golf-1} becomes
$$
v_2{y_2'}^2z_0^4   -x_4z_1^2+ v_2{y_3'}^2=0.
$$
At this point it is convenient to deduce a coprimality
condition which follows from the assumption made at the outset that
$\hcf(x_0,\ldots,x_4)=1$.
Recalling the various changes of variables
that we have made so far, we easily conclude that
\begin{equation}\lab{polo-2}
\hcf(v_2y_2',x_4)=1.
\end{equation}
But then it follows that $v_2$ must divide $z_1^2$ in the above
equation.  Since $v_2$ is square-free and positive we may conclude that there
exists $y_1' \in \N$ such that
$$
z_1=v_2y_1'.
$$
This leads to the equation
\begin{equation}\lab{golf-2}
{y_2'}^2z_0^4   -v_2x_4{y_1'}^2+ {y_3'}^2=0,
\end{equation}
in addition to the further coprimality condition
\begin{equation}\lab{polo-3}
\hcf(z_0,v_2y_1')=1,
\end{equation}
that follows from \eqref{polo-1}.

Next we let $v_1=\hcf(y_1',y_3')$.  Then
$v_1 \in \N$ and we may write
\begin{equation}\lab{subs-1}
z_0=y_0, \quad y_1'=v_1y_1, \quad y_3'=v_1y_3,
\end{equation}
for $y_0,y_1,y_3 \in \N$.
Under this change of variables we deduce
from \eqref{polo-3} and the definitions of $v_1,v_2$, that
\begin{equation}\lab{polo-4'}
|\mu(v_2)|=1, \quad \hcf(y_1,y_3)=\hcf(y_0,v_1v_2y_1)=1,
\end{equation}
where $\mu(n)$ denotes the M\"obius function for any $n \in \N$.
We proceed by substituting \eqref{subs-1} into \eqref{golf-2}.  This
leads to the equation
$$
y_0^4{y_2'}^2   -v_1^2v_2x_4y_1^2+ v_1^2y_3^2=0,
$$
from which it follows that $v_1^2
\mid y_0^4{y_2'}^2$.  In view of \eqref{polo-4'}
it follows that there exist~$y_2 ,y_4 \in \N$ such
that
$$
x_4=y_4, \quad y_2'=v_1y_2.
$$
Our investigation has therefore produced the equation
\begin{equation}\lab{ut}
y_0^4y_2^2   -v_2y_1^2y_4+ y_3^2=0,
\end{equation}
and \eqref{polo-2} becomes
\begin{equation}\lab{polo-5}
\hcf(y_4,v_1v_2y_2)=1.
\end{equation}

We take a moment to deduce two further coprimality conditions
\begin{equation}\lab{got}
\hcf(y_2,y_3)=1,\qquad\hcf(y_2,v_2y_1)=1.
\end{equation}
      Suppose
that there exists a prime divisor $p$ common to both $y_2$ and
$y_3$.  Then it follows from~\eqref{ut} that $p^2 \mid v_2y_1^2y_4$.
But \eqref{polo-4'} and \eqref{polo-5} together imply that $\hcf(p,
y_1y_4)=1$.  Hence $p^2 \mid v_2$, which is impossible since $v_2$ is
square-free.
The equation \eqref{ut} and the fact that
$\hcf(y_2,y_3)=1$ clearly yield the second relation in \eqref{got}.
Combining \eqref{got} with \eqref{polo-4'} and
\eqref{polo-5},  we therefore obtain the relations
\begin{equation}\lab{polo-6}
\hcf(y_0,v_1v_2y_1)=
\hcf(y_3,y_1y_2)=\hcf(y_4,v_1v_2y_2)=1,
\end{equation}
and
\begin{equation}\lab{polo-7}
|\mu(v_2)|=1, \quad \hcf(y_2,v_2y_1)=1.
\end{equation}

At this point we may summarise our argument as follows.
Let $\mcal{T}$ denote the set of
$
(\ma{v},\ma{y} )=(v_1,v_2,  y_0, \ldots,y_4) \in \N^7
$
such that \eqref{ut}, \eqref{polo-6} and \eqref{polo-7}
hold.  Then for any $\x \in N^5$ counted by $N(Q_1,Q_2;B)$, we have
shown that there exists $(\ma{v},\ma{y} ) \in \mcal{T}$ such that
$$
\begin{array}{l}
x_0 = v_1^2v_2y_0^2y_2^2, \quad
x_1 = v_1^4v_2^3 y_1^2  y_2^2, \quad
x_2 = v_1^3v_2^2y_0y_1y_2^2, \\
x_3 = v_1^2v_2  y_2 y_3, \quad x_4=y_4.
\end{array}
$$
Conversely, given any 
$(\ma{v},\ma{y} ) \in \mcal{T} $, the point $\x$ given above
will be a solution of the equations $Q_1(\x)=Q_2(\x)=0$, with
$\x \in N^5$.  To check that $\hcf(x_0,\ldots,x_4)=1$,
we first observe that
$$
\hcf(x_0,x_1,x_2)=v_1^2v_2y_2^2
\hcf(y_0^2, v_1^2v_2^2 y_1^2, v_1v_2y_0y_1)=v_1^2v_2y_2^2,
$$
by \eqref{polo-6}.  But then \eqref{polo-6} also implies that
$$
\hcf(x_0,x_1,x_2,x_3)=v_1^2v_2 y_2 \hcf(y_2,y_3)=v_1^2v_2 y_2,
$$
whence
$$
\hcf(x_0,\ldots,x_4)=\hcf(v_1^2v_2y_2,y_4)=1.
$$
On observing that the inequality
$y_4\leq B$ is equivalent to  $y_0^4y_2^2+y_3^2\leq
Bv_2y_1^2$, by~\eqref{ut}, we therefore conclude the proof of the
following result.

\begin{lem}\lab{base}
Let $B\geq 1$.  Then we have
$$
N(Q_1,Q_2;B)= \#\Big\{(\ma{v,y}) \in \mcal{T}: v_1^4v_2^3 y_1^2
       y_2^2 \leq B, ~y_0^4y_2^2+y_3^2\leq Bv_2y_1^2\Big\}.
$$
\end{lem}

It will become clear in subsequent sections that the equation
\eqref{ut} is a pivotal ingredient in our proof of Theorems
\ref{main'} and \ref{main}.

\section{The main offensive}\lab{final}

In this section we use Lemma \ref{Reduc1} and Lemma \ref{base} to
obtain an initial estimate for~$N_{U,H}(B)$, which will then be used
to deduce the statement of Theorem~\ref{main'} in~\S\ref{height}.
Before beginning this task, it will be helpful to first outline
our strategy.  It follows from the statement of Lemma \ref{base} that
any vector $(\ma{v},\ma{y})$ counted by $N(Q_1,Q_2;B)$ is 
constrained to lie in the region defined by \eqref{ut} and the
inequalities
$$
v_1^4v_2^3 y_1^2y_2^2 \leq B, \quad y_0^4y_2^2+y_3^2\leq Bv_2y_1^2.
$$

The bulk of our work will be taken up with handling the summation over
the variables $y_0,y_3$ and $y_4$, for fixed values of $(\ma{v},y_1,y_2)$.
The essential trick will be to view \eqref{ut} as a congruence
condition in order to handle the summation over $y_4$. Then we will
need to count values of $y_3$ such that 
\begin{equation}\lab{Y3}
y_3 \leq \sqrt{Bv_2y_1^2-y_0^4y_2^2}=Y_3, 
\end{equation}
say, subject to certain coprimality conditions, and then
finally values of $y_0$ such that
\begin{equation}\lab{Y0}
y_0 \leq \frac{B^{1/4}v_2^{1/4}y_1^{1/2}}{y_2^{1/2}}=Y_0, 
\end{equation}
say, subject to certain coprimality conditions.
The overall contribution from $y_0,y_3$ and $y_4$, which we henceforth denote by
$S(\ma{v},y_1,y_2)$, will be estimated in \S \ref{suby3y4}.
Finally in \S \ref{sublast} we will sum $S(\ma{v},y_1,y_2)$ over
the remaining values of $\ma{v},y_1,y_2$ subject to the inequality
\begin{equation}
v_1^4v_2^3y_1^2y_2^2\leq B,
\lab{ineg-3}
\end{equation}
and certain coprimality conditions.

\subsection{Congruences and equidistribution}

It will be convenient to collect together some technical facts about
congruences that will be needed in \S \ref{suby3y4} and 
\S\ref{sublast}. We begin by
discussing the arithmetic function $\eta(q)$,
defined to be the number of square roots of $-1$ modulo $q$.
The Chinese Remainder Theorem implies that $\eta(q)$ is
multiplicative, and for any $\nu \geq 1$ we have
$$
\eta(p^\nu)=
\left\{
\begin{array}{ll}
2  & \hbox{if } p\equiv 1  \mod{4},\\
0  & \hbox{if } p\equiv 3  \mod{4},\\
1  & \hbox{if } p=2, \,\nu=1,\\
0  & \hbox{if } p=2, \,\nu \geq 2.
\end{array}
\right.
$$
Let $\omega(n)$ denote the number of distinct prime factors of $n \in
\N$, and recall the definition of $\chi$, the real non-principal
character modulo $4$.
Then it is not difficult to see that we have
\begin{equation}\lab{triv}
\eta(q) \leq  \sum_{d \mid q}|\mu(d)| \chi(d) \leq 2^{\omega(q)},
\end{equation}
for any positive integer $q$.

Define the real-valued function $\psi(t)=\{t\}-1/2$, where
$\{\al\}$ denotes the fractional part of $\al \in \R$.  Then
$\psi$ is periodic with period $1$.  We proceed by recording the
following simple estimate.

\begin{lem}\lab{cong}
Let $a, q \in \Z$ be such that $q>0$, and let $t\geq 0$.
Then we have
$$
\#\{0<n \leq t: n\equiv a\mod{q} \}=\frac{t}{q}+ r(t;a,q),
$$
where
$$
r(t;a,q)=\psi\Big( -\frac{a}{q}\Big) - \psi\Big( \frac{t-a}{q}\Big).
$$
\end{lem}

\begin{proof}
This follows on taking $t_1=0$ and $t_2=t$ in \cite[Lemma 3]{1}.
\end{proof}

We will also need to prove a result about the average order of the
function $\psi$, that plays the same role in this work that 
\cite[Lemma 4]{1} did in our previous work.  
A crucial ingredient in this will be the following Diophantine
approximation result, which may be of independent interest.

\begin{lem}\lab{approxdio}
Let $ q\in \N$ and let $\vr$ be a square root of $-1$ modulo $q$.  For each 
non-zero integer $b$, there exist
  coprime integers $u,v$  such that
\begin{equation}\lab{rat}
\Big| \frac{b \vr}{q} - \frac{u}{v}\Big| \leq \frac{1}{v\sqrt{2q}}
\end{equation} and
\begin{equation}\lab{ineg-v}
\frac{\sqrt{q/2}}{|b|} \leq v \leq \sqrt{2q}.
\end{equation}
\end{lem}

\begin{proof}
By Dirichlet's approximation theorem we may find coprime integers $u,v$ such that $1\leq v \leq 
\sqrt{2q}$ and \eqref{rat} holds. We claim that any such $v$
automatically satisfies \eqref{ineg-v}, for which it clearly remains to
establish the lower bound.  To do so we first note 
that $|r| \leq \sqrt{q/2}$, where $r$ denotes the residue of $b \vr v$ modulo
$q$.  But then it follows that
$$
r^2 \equiv b^2\vr^2v^2 \equiv -b^2 v^2 \mod{q},
$$
since $\vr$ is a square-root of $-1$ modulo $q$, whence $q \mid
r^2+b^2v^2$.  Thus we obtain the system of inequalities
$$
q \leq r^2+b^2v^2 \leq q/2 + b^2v^2,
$$
which thereby gives the lower bound in \eqref{ineg-v}.
\end{proof}

Given any real-valued function $f$ defined on an interval $I$, and
given coefficients $c_1,c_2\in \R$ and $q\in\N$, define the sum
$$
S_I(f;c_1,c_2,q)=
\sum_{x\in \Z\cap I}
\psi\Big(\frac{c_1 f(x)-c_2 x^2}{q}\Big).
$$
We proceed by establishing the following result, which permits
$S_I(f;c_1,c_2,q)$ to be estimated over short intervals $I$ under suitable hypotheses.

\begin{lem}\lab{av-order}
Let $b,q \in \Z$ with $q>0$ and $b \neq 0$, and let
$\vr$ be a square root of $-1$ modulo $q$.  Let $\gamma \in \R$. 
Given a bounded interval
$I\subset \R$, let $f:I\rightarrow \R$
be a continuously differentiable function on $I$. Write
$$
\la(f)=\int_I \big|f'(t) \big| \d t  ,
$$
and define
$$
g_{\al,\be}(t)=f(t)-\al-\be t,
$$
for any $\al,\be \in \R$.  Then for any $H\geq 1$ we have
\begin{align*}
S_I(f;\gamma,b\vr,q)
\ll 
\frac{\m(I)}{H}+\log (H R) \Big(1+\frac{|\gamma|H\la(g_{\al,\beta})}{q}  \Big)  
\Big(q^{1/4}+\frac{|b|^{1/2}H^{1/2}\m(I)}{q^{1/4}}\Big),
\end{align*}
where $\m(I)=\meas(I)+1$ and $R$ is such that $I\subseteq [-R,R]$.
\end{lem}

\begin{proof}
We employ Vaaler's trigonometric polynomial approximation
\cite{vaaler} to $\psi$.  This implies that for any $H \geq 1$, there exist
coefficients $c_h$ such that
$$
\psi(t) \leq  \sum_{0<|h|\leq H} c_h e(ht)+  O\Big(\frac{1}{H}\Big),
$$
with $c_h \ll 1/|h|$ and $e(\theta)=\exp(2\pi i \theta)$. 
It therefore  follows that
\begin{equation}\lab{step1}
S_I(f;\gamma,b\vr,q) \ll  \frac{\m(I)}{H}+
\sum_{0<|h|\leq H} \frac{1}{|h|} \Big|\sum_{x\in \Z\cap I}  
e\Big(\frac{h\gamma f(x)-bh\vr x^2}{q}\Big)
\Big|,
\end{equation}
for any $H\geq 1$.

Define the sum 
$$
S_{h,Y}(F)= \sum_{ Y/2 <y\leq Y} 
e\Big(\frac{hF(y)-bh\vr y^2}{q}\Big),
$$
for fixed values of $h \in \N$, $Y \geq 1$ and any real valued function
$F:(Y/2,Y]\rightarrow \R$.
We begin by estimating the sum when $F(t)=\al+\be t$. 
Mimicking the proof of Weyl's
inequality, we obtain
\begin{align*}
|S_{h,Y}(\al+\be t)|^2 &= \Big| \sum_{\colt{Y/2 <y\leq Y}{Y/2 <y+z\leq
  Y}} e\Big(\frac{-\beta h z+bh\vr (z^2+2yz)}{q}\Big)\Big|\\
&\leq \sum_{|z|\leq Y}
\Big| \sum_{y\in I(z)} 
e\Big(\frac{2bh\vr yz}{q}\Big)\Big|,
\end{align*}
for some interval $I(z)$ of length $O(Y)$, depending on $z$.
For any $\al\in \R$, let $\|\al\|$ denote the
distance from $\al$ to the nearest integer. Then it follows that
\begin{align*}
|S_{h,Y}(\al+\be t)|^2
&\ll \sum_{|z|\leq Y}
\min\Big\{ Y, \frac{1}{\|2b \vr hz/q\|}\Big\}.
\end{align*}

For any $\al \in \R$ and any $\delta>0$, let $T(N;\al,\delta)$
denote the number of positive integers $n \leq N$ such that $\|\al n\|
\leq \delta$.  A result due to Heath-Brown 
\cite[Lemma 6]{weyl} yields
\begin{equation}\lab{hb}
T(N;\al,\delta) \leq 4(1+v \delta)\Big(1+\frac{N}{v}\Big),
\end{equation}
provided that $\al$ has a rational approximation $|\al-u/v| \leq
1/v^{2}$ with coprime integers $u,v$ such that $v \neq 0$.
On combining \eqref{hb} with Lemma \ref{approxdio}, we 
therefore deduce that
\begin{align*}
|S_{h,Y}(\al+\be t)|^2
&\ll
\max_{1/Y  \leq \delta \leq 1}\frac{ T(Y;2bh\vr/q, \delta)}{\delta}\\
&\ll (Y+v)\Big(1+\frac{Y}{v} \Big)\\
&\ll \Big(q^{1/4} + \frac{|bh|^{1/2}Y}{q^{1/4}}\Big)^2,
\end{align*}
whence
$$
S_{h,Y}(\al+\be t)
\ll
q^{1/4} + \frac{|bh|^{1/2}Y}{q^{1/4}}.
$$

Given any $\gamma \in \R$ and any function $f: I\rightarrow  \R$ as in the statement of the lemma,
it now follows from an application of partial summation that 
\begin{align*}
S_{h,Y}(\gamma f)&= S_{h,Y}\big(\gamma g_{\al,\be}(t) + \alpha'+\beta' t\big)\\
&\ll
\Big(1+ \frac{h|\gamma| \lambda(g_{\alpha, \beta})}{q}\Big)
\Big(q^{1/4} + \frac{|bh|^{1/2}Y}{q^{1/4}}\Big),
\end{align*}
with $\alpha'=\alpha \gamma$ and $\beta'=\beta \gamma$.
Finally we employ this estimate in \eqref{step1}, and sum over dyadic
intervals for $Y\ll R$, in order to complete the proof
of Lemma \ref{av-order}.
\end{proof}

\subsection{Summation over the variables $y_0,y_3$ and $y_4$}\lab{suby3y4}

Let $(\ma{v},y_1,y_2) \in \N^4$ satisfy
\eqref{polo-7} and \eqref{ineg-3}.
As indicated above, we will denote the triple summation over $y_0,y_3$
and $y_4$ by $S(\ma{v},y_1,y_2)$.
In this summation we will need to take into account the coprimality conditions
\eqref{polo-6}.  We begin by treating the condition
$\hcf(y_4, v_1v_2)=1$.  Note here that the condition $\hcf(y_4,y_2)=1$
follows immediately from \eqref{ut} and the remaining conditions. A
M\"obius inversion yields
$$
S(\ma{v},y_1,y_2)=\sum_{ k_4 \mid v_1  v_2} \mu (k_4)S_{k_4},
$$
where the definition of $S_{k_4}$ is as for $S(\ma{v},y_1,y_2)$ but
with the extra
condition that $k_4\mid y_4$ and without the
condition $\hcf(y_4, v_1v_2y_2)=1$.
Now it  straightforward to deduce from \eqref{ut}, \eqref{polo-6} and \eqref{polo-7},
that
$$
\hcf(k_4,y_2) = \hcf(v_1v_2,y_2,y_3^2) =1,
$$
for any $k_4$ dividing $v_1v_2$ and $y_4$.
It follows that
$$
S(\ma{v},y_1,y_2)=\sum_{\colt{y_0\leq
Y_0}{\hcf(y_0,v_1v_2y_1)=1}}\sum_{\colt{k_4 \mid
v_1v_2}{\hcf (k_4,y_2)=1}}\mu(k_4)
S_{k_4}',
$$
where
$$
S_{k_4}'=
\#\left\{ y_3\in\N:
\hcf(y_3,y_1y_2)=1,~y_3\leq Y_3,~
y_3^2\equiv -y_0^4y_2^2 \mod {k_4v_2y_1^2}
\right\},
$$
and $Y_3$ is given by \eqref{Y3}. Let $\vr$ be a solution of
the congruence 
$$
\vr^2\equiv -1 \pmod{k_4v_2y_1^2}.
$$ 
Then a little
thought reveals that 
\begin{equation}\lab{ham}
\hcf(\vr y_0^2y_2,k_4v_2y_1^2)=1,
\end{equation}
for any $k_4 \mid v_1v_2$ such that $\hcf(k_4,y_2)=1$.  We may
therefore conclude that 
$$
S(\ma{v},y_1,y_2)=\sum_{\colt{y_0\leq
Y_0}{\hcf(y_0,v_1v_2y_1)=1}}\sum_{\colt{k_4 \mid
v_1v_2}{\hcf (k_4,y_2)=1}}\mu(k_4)
\sum_{\colt{\vr \leq k_4v_2y_1^2}{\vr^2\equiv -1\mod
{k_4v_2y_1^2}}}S_{k_4}(\vr),
$$
where
$$
S_{k_4}(\vr)=
\#\left\{ y_3\in\N:
\hcf(y_3,y_2)=1,~y_3\leq Y_3,~
y_3\equiv \vr y_0^2y_2 \mod {k_4v_2y_1^2}
\right\}.
$$

In order to estimate $S_{k_4}(\vr)$ we employ a further M\"obius
inversion to treat the coprimality condition $\hcf(y_3,y_2)=1$.
Thus we have
$$
S_{k_4}(\vr)=
\sum_{k_2\mid y_2 } \mu (k_2)
S_{k_2,k_4 }(\vr),
$$
where
$$
S_{k_2,k_4 }(\vr)=
\#\left\{ y'_3\in\N :
y'_3\leq Y_3/k_2, \, k_2y'_3\equiv \vr y_0^2y_2 \mod {k_4v_2y_1^2}
\right\}.
$$
We must therefore estimate the number of positive integers
contained in a certain interval, which belong to a certain arithmetic
progression.
For fixed values of $\ma{v},y_0,y_1,y_2,k_2,k_4,\vr$, it follows from
\eqref{ham} that there exists a unique positive integer
$b \leq k_4v_2y_1^2$ such that $\hcf(b,k_4v_2y_1^2)=1$ and
\begin{equation}\lab{*}
bk_2\equiv \vr y_2\mod{k_4v_2y_1^2}.
\end{equation}
We may therefore employ Lemma \ref{cong} to deduce that
\begin{align*}
S_{k_2,k_4 }(\vr)&=
\#\left\{ y'_3\in\N :
y'_3\leq Y_3/k_2, \, y'_3\equiv by_0^2 \mod {k_4v_2y_1^2}
\right\}\\
&= \frac{Y_3}{k_2k_4v_2y_1^2}+r(Y_3/k_2;
by_0^2,k_4v_2y_1^2),
\end{align*}
where
\begin{equation}\lab{viva}
r(Y_3/k_2; by_0^2,k_4v_2y_1^2)=\psi\Big( \frac{-by_0^2}{k_4v_2y_1^2}\Big)
-\psi\Big(\frac{Y_3/k_2-by_0^2}{k_4v_2y_1^2}\Big).
\end{equation}

Putting everything together we have therefore established that
\begin{equation}\lab{manray}
S(\ma{v},y_1,y_2)= \vt(\ma{v},y_1,y_2)\hspace{-0.2cm}
\sum_{\colt{y_0\leq Y_0}{\hcf(y_0,v_1v_2y_1)=1}}
\hspace{-0.2cm}
\frac{Y_3}{v_2y_1^2}+ R_1(B;\ma{v},y_1,y_2),
\end{equation}
with
\begin{align*}
\vt(\ma{v},y_1,y_2)
&=\sum_{\colt{k_4\mid v_1 v_2 }{ \hcf (k_4,y_2)=1}}  \frac{\mu
(k_4)\eta(k_4v_2y_1^2)}{k_4} \sum_{k_2\mid y_2 }\frac{\mu (k_2)}{k_2}\\
&=\frac{\phi(y_2)}{y_2}\sum_{\colt{k_4\mid v_1 v_2 }{ \hcf
(k_4,y_2)=1}}  \frac{\mu
(k_4)\eta(k_4v_2y_1^2)}{k_4},
\end{align*}
and
\begin{equation}\lab{R_1}
R_1(B;\ma{v},y_1,y_2)=
\hspace{-0.2cm}
\sum_{\colt{k_4\mid v_1 v_2 }{ \hcf (k_4,y_2)=1}}
\hspace{-0.2cm}
\mu(k_4)
\hspace{-0.2cm}
\sum_{\colt{\vr \leq k_4v_2y_1^2}{\vr^2\equiv -1\mod {k_4v_2y_1^2}}}
\sum_{k_2\mid y_2}
\mu (k_2)
F_1(B).
\end{equation}
Here we have set
$$
F_1(B)=F_1(B;\ma{v},y_1,y_2,k_2,k_4,\vr)=\sum_{\colt{y_0\leq
Y_0}{\hcf(y_0,v_1v_2y_1)=1}}
r(Y_3/k_2; by_0^2,k_4v_2y_1^2),
$$
and a straightforward calculation reveals that
\begin{equation}\lab{calcultheta}
\vt(\ma{v},y_1,y_2)=
\frac{\eta(v_2y_1^2)\phi(y_2)}{ y_2} \!\!\!
\prod_{\colt{p\mid v_1}{p\nmid v_2y_1y_2}}
\!\!\!\!\!\Big(1-\frac{1+\chi(p)}{p}\Big)\!\!\!\!\!
\prod_{ p\mid v_2\hcf(v_1,y_1) }
\!\!\!\!\!\Big(1-\frac{\chi(p)}{p}\Big).
\end{equation}
We have used here the fact that whenever
$\eta(v_2y_1^2) \neq 0$, any prime divisor of $v_2\hcf(v_1,y_1)$ must
be congruent to $1$ modulo $4$.
We proceed by establishing the following result.

\begin{lem}\lab{sumR1}
Let $\ve>0$.  Then for any $B\geq 1$ we have
$$
\sum_{\colt{v_1,v_2,y_1,y_2\in 
\N}{\mbox{\scriptsize{\eqref{polo-7},\eqref{ineg-3}
hold}}}} R_1(B;v_1,v_2,y_1,y_2)
\ll_\ve B^{17/20+\ve}.
$$
\end{lem}

\begin{proof}
We begin with a M\"obius inversion to remove the coprimality condition
in the summation over $y_0$ in the definition of $F_1(B)$.
Thus it follows from \eqref{viva} that
$F_1(B)=F_1^{\mathrm{a}}(B)-F_1^{\mathrm{b}}(B)$, with 
$$
F_1^{\mathrm{a}}(B)= \sum_{k_0\mid v_1v_2y_1} \mu(k_0) 
\sum_{y_0\leq Y_0/k_0} \psi\Big( \frac{-bk_0^2y_0^2}{k_4v_2y_1^2}\Big),
$$
and 
$$
F_1^{\mathrm{b}}(B)= \sum_{k_0\mid v_1v_2y_1} \mu(k_0) 
\sum_{y_0\leq Y_0/k_0} \psi\Big(\frac{Y_3/k_2-bk_0^2y_0^2}{k_4v_2y_1^2}\Big).
$$
On writing $y_2=k_2y_2'$, \eqref{*} implies that
$
b \equiv \vr y_2' ~(\bmod~{k_4v_2y_1^2}).
$
Here we recall that $\vr$ is a square root of $-1$ modulo $k_4v_2y_1^2$.
We now set
\begin{equation}
  \label{eq:qq}
q=k_4v_2{y_1}^2.
\end{equation}
An application of Lemma \ref{av-order} with $\gamma=0$ 
yields
\begin{align*}
F_1^{\mathrm{a}}(B)
&\ll_\ve \sum_{k_0 \mid v_1v_2y_1} B^\ve 
  \Big( \frac{Y_0}{H}+ q^{1/4} + \frac{{y_2'}^{1/2}H^{1/2}Y_0}{q^{1/4}} \Big)\\
&\ll_\ve
B^\ve  \Big(\frac{Y_0}{H}+(k_4v_2 y_1^2)^{1/4}+\frac{y_2^{1/2}H^{1/2}Y_0}{(v_2 y_1^2)^{1/4}}\Big),
\end{align*}
for any $H \geq 1$. 
Substituting this into \eqref{R_1}, we therefore obtain an overall
contribution of 
\begin{align*}
&\ll_\ve B^\ve
\Big(\frac{Y_0}{H}+v_1^{1/4} v_2^{1/2}y_1^{1/2}+\frac{y_2^{1/2}H^{1/2}Y_0}{v_2^{1/4} y_1^{1/2}}
\Big).
\end{align*}
to $R_1(B;\ma{v},y_1,y_2)$ from $F_1^{\mathrm{a}}(B)$.
Here, we have used the estimate \eqref{triv} for $\eta(k_4v_2y_1^2)$.
It remains to sum this estimate over all $v_1,v_2,y_1,y_2 
\in\N$ satisfying \eqref{polo-7} and \eqref{ineg-3}.
Recall the definition \eqref{Y0} of $Y_0$.
Then for any $\ve>0$ we obtain the overall contribution 
\begin{align*}
&\ll_\ve B^\ve 
\sum_{\colt{v_1,v_2,y_1,y_2\in \N}{y_1 \leq
B^{1/2}/(v_1^2v_2^{3/2}y_2)}} 
\Big(\frac{B^{1/4}v_2^{1/4}y_1^{1/2}}{Hy_2^{1/2}}+v_1^{1/4} v_2^{1/2}y_1^{1/2}
+B^{1/4}H^{1/2}\Big)\\
&\ll_\ve B^\ve \Big(\frac{B}{H}+B^{3/4}H^{1/2}\Big)
\end{align*}
to $\sum_{v_1,v_2,y_1,y_2} R_1(B;\ma{v},y_1,y_2)$ from the term
$F_1^{\mathrm{a}}(B)$. This term therefore contributes $O_\ve(B^{5/6+\ve})$ to Lemma
\ref{sumR1}, on taking $H=B^{1/6}$, which is satisfactory.

We must now turn to the overall contribution from the term
$F_1^{\mathrm{b}}(B)$.  Define the function
$f(t)= \sqrt{1-t^4}$. Then it follows from \eqref{Y3} that 
$$
\frac{Y_3}{k_2}=\sqrt{Bv_2y_1^2}f(y_0/Y_0'),
$$
in the above expression for $F_1^{\mathrm{b}}(B)$, where we have set
$Y_0'=Y_0/k_0$ for short.
We will estimate $F_1^{\mathrm{b}}(B)$ by breaking the summation over
$y_0$ into intervals of length $O(Y_0'/M)$, for any $M\geq 1$ to be
selected in due course. Let 
$$
U_m=
\sum_{mY_0'/M<y_0\leq (m+1)Y_0'/M}
\psi\Big(\frac{\sqrt{Bv_2y_1^2}f(y_0/Y_0')-\vr y_2'k_0^2y_0^2}{q}\Big),
$$
for each $m \in \N$, where $q$ is given by \eqref{eq:qq}.
We wish to apply 
Lemma \ref{av-order} to estimate $U_m$, in which we will take 
$\gamma=\sqrt{Bv_2y_1^2}$.

Let us begin by handling $U_m$ for $m \leq M-2$. Here we take 
$$
\beta=\frac{f'(m/M)}{M}, \quad \alpha=f(m/M)-\frac{mf'(m/M)}{M^2},
$$
and proceed to consider the size of $g_{\al,\be}(t/Y_0')$ on the range
of summation.
An application of the mean value theorem implies that for $mY_0'/M<t\leq
(m+1)Y_0'/M$ we have 
\begin{align*}
g_{\al,\be}(t/Y_0')&= f(t/Y_0')-f(m/M)-\frac{1}{M}f'(m/M)
\Big(\frac{t}{Y_0'}-\frac{m}{M}\Big)\\
& \ll \frac{1}{M^2}\sup_{\tau\in (m,m+1]}|f''(\tau/M)| \\
&\ll\frac{(1-(m+1)^4/M^4)^{-3/2}}{M^2} .
\end{align*}
On making the change of variables $m=M-m'-1$, so that 
$$
(1-(m+1)^4/M^4)^{-3/2}\ll (M/m')^{3/2}, 
$$
we therefore deduce that
\begin{align*}
\sum_{1\leq m\leq M-2} \lambda(g_{\al,\be}(t/Y_0')) 
&\ll  M^{-2}\sum_{1\leq m\leq M-2}(1-(m+1)^4/M^4)^{-3/2}\\
&\ll  M^{-1/2}\sum_{1\leq m'\leq M-2}{m'}^{-3/2} \ll M^{-1/2}.
\end{align*} 
It now follows from Lemma \ref{av-order} that
\begin{align*}
\sum_{1\leq m\leq M-2}U_m 
&\ll_\ve B^\ve
\Big(\frac{Y_0}{H}+\Big(M+\frac{B^{1/2}H}{v_2^{1/2}y_1M^{1/2}}\Big)\Big(v_1^{1/4}
v_2^{1/2}y_1^{1/2}+\frac{y_2^{1/2}H^{1/2}Y_0}{v_2^{1/4} y_1^{1/2}M}\Big)
\Big).
\end{align*}

To estimate $U_m$ for $m=M-1$ or $M$, we choose
$\beta= \alpha=0$ in Lemma \ref{av-order}. In particular it follows
that 
$g_{\al,\be}(t/Y_0')= f(t/Y_0')\ll M^{-1/2}$ and so
$$
\lambda(g_{\al,\be}(t/Y_0')) \ll  M^{-1/2},
$$
for $mY_0'/M<t\leq (m+1)Y_0'/M$ and $M-2<m \leq M $.
We therefore deduce that 
\begin{align*}
\sum_{M-2< m\leq M}U_m 
&\ll_\ve B^\ve
\Big(\frac{Y_0}{H}+\Big(1+\frac{B^{1/2}H}{v_2^{1/2}y_1M^{1/2}}\Big)\Big(v_1^{1/4}
v_2^{1/2}y_1^{1/2}+\frac{y_2^{1/2}H^{1/2}Y_0}{v_2^{1/4} y_1^{1/2}M}\Big)
\Big),
\end{align*}
which once combined with the above leads to the conclusion that
\begin{align*}
 F_1^{\mathrm{b}}(B) 
&\ll_\ve B^\ve
\Big(\frac{Y_0}{H}+\Big(M+\frac{B^{1/2}H}{v_2^{1/2}y_1M^{1/2}}\Big)\Big(v_1^{1/4}
v_2^{1/2}y_1^{1/2}+\frac{y_2^{1/2}H^{1/2}Y_0}{v_2^{1/4} y_1^{1/2}M}\Big)
\Big)\\
&\ll_\ve B^\ve\Big(
 \frac{Y_0}{H}+ \frac{B^{1/3}H^{2/3}}{v_2^{1/3}y_1^{2/3} } \Big(v_1^{1/4}
v_2^{1/2}y_1^{1/2}+\frac{y_2^{1/2}  Y_0 v_2^{1/12}y_1^{1/6}}{ 
B^{1/3}H^{1/6}}
\Big)
\Big),
\end{align*}
on choosing $M=H^{2/3}B^{1/3}/(v_2^{1/3}y_1^{2/3})$. Note that $M\geq
1$ since $v_2 y_1^2\leq B$.  
Recall the definition \eqref{Y0} of $Y_0$ and let $\ve>0$. Then
arguing as above, we conclude that we have an overall contribution 
\begin{align*}
\ll_\ve B^\ve 
\sum_{v_1,v_2,y_1,y_2} 
\Big(\frac{B^{1/4}v_2^{1/4}y_1^{1/2}}{Hy_2^{1/2}}
+\frac{v_1^{1/4}v_2^{1/6}
B^{1/3}H^{2/3}}{y_1^{1/6}} + {B^{1/4}H^{1/2}}  \Big)
\end{align*}
to $\sum_{v_1,v_2,y_1,y_2} R_1(B;\ma{v},y_1,y_2)$ from the term
$F_1^{\mathrm{b}}(B)$. Here the summation is over 
$v_1,v_2,y_1,y_2\in\N$ such that $y_1 \leq
B^{1/2}/(v_1^2v_2^{3/2}y_2)$.
Taking $H=y_1^{3/10}$, we therefore obtain the satisfactory contribution
\begin{align*}
&\ll_\ve B^\ve 
\sum_{v_1,v_2,y_1,y_2} 
\Big(\frac{B^{1/4}v_2^{1/4}y_1^{1/5}}{y_2^{1/2}}
+v_1^{1/4}v_2^{1/6}y_1^{1/30}
B^{1/3} + {B^{1/4}y_1^{3/20}}  \Big)\\
&\ll_\ve B^\ve 
\sum_{v_1,v_2,y_1,y_2} 
\Big(\frac{B^{1/4}v_2^{1/4}y_1^{1/5}}{y_2^{1/2}}
+v_1^{1/4}v_2^{1/6}y_1^{1/30}
B^{1/3} \Big)\\
&\ll_\ve B^{17/20+\ve},
\end{align*}
on noting that ${B^{1/4}y_1^{3/20}} \leq B^{1/4+3/40} \leq B^{1/3}$
and summing first over $y_1$. 
This completes the proof of Lemma \ref{sumR1}.
\end{proof}

Recall the definitions \eqref{Y3}, \eqref{Y0} of $Y_0, Y_3$, and write
$f(t)= \sqrt{1-t^4}$ as above. Our next task is to examine the sum
$$
\sum_{\colt{y_0\leq Y_0}{\hcf(y_0,v_1v_2y_1)=1}}
\frac{Y_3}{v_2y_1^2} =  \frac{B^{1/2}}{v_2^{1/2}y_1}\sum_{\colt{y_0\leq
Y_0}{\hcf(y_0,v_1v_2y_1)=1}}
f( y_0 /Y_0),
$$
that appears in \eqref{manray}.  On performing a M\"obius inversion
we easily deduce that
\begin{align*}
\sum_{\colt{y_0\leq
Y_0}{\hcf(y_0,v_1v_2y_1)=1}}f( y_0 /Y_0)&=
\sum_{k_0 \mid v_1v_2y_1}\mu(k_0)
\sum_{y_0'\leq Y_0'}  f( y_0' /Y_0'),
\end{align*}
where $Y_0'=Y_0/k_0$. An application of partial summation, together
with the change of variables $t=u^{1/4}Y_0'$, reveals  that
\begin{align*}
\sum_{y_0'\leq Y_0'}  f(y_0'/Y_0') &=[Y_0']f(1)-\int_0^{Y_0'}
[t]f'\Big(\frac{t}{Y_0'}\Big)
\frac{\d t}{Y_0'}\\ &=- \int_0^{1}
\big[u^{1/4}Y_0'\big]f'(u^{1/4})\frac{\d u}{4u^{3/4}}  =c
Y_0'+F_2(Y_0'),
\end{align*}
where
$$
F_2(Y_0')=-\frac{1}{2}\int_0^{1}\big\{ u^{1/4}Y_0'  \big\}\frac{\d
u}{\sqrt{1-u}},
$$
and
\begin{equation}\lab{TAU}
c=\int_0^1 \frac{u^{1/4}\d u}{ 2\sqrt{1-u}}.
\end{equation}
We easily conclude that
\begin{align*}
\sum_{\colt{y_0\leq Y_0}{\hcf(y_0,v_1v_2y_1)=1}}
\frac{Y_3}{v_2y_1^2}
&=  c\frac{\phi(v_1v_2y_1)}{v_1v_2y_1}
\frac{B^{3/4}}{v_2^{1/4}y_1^{1/2}y_2^{1/2}}
+\frac{R_2(B;\ma{v},y_1,y_2)}{\vt(\ma{v},y_1,y_2)},
\end{align*}
where
\begin{equation}\lab{R_2}
R_2(B;\ma{v},y_1,y_2)= \vt(\ma{v},y_1,y_2)
\frac{B^{1/2}}{v_2^{1/2}y_1}\sum_{k_0\mid v_1v_2y_1}\mu(k_0)F_2(Y_0').
\end{equation}
We now need to sum $R_2(B;\ma{v},y_1,y_2)$ over all of the relevant values
of $y_2$.  Let
$$
Y_2=\frac{B^{1/2}}{v_1^2v_2^{3/2}y_1}.
$$
Then the following result holds.

\begin{lem}\lab{sumR2}
Let $\ve>0$ and let $B\geq 1$.  Then we have
$$
\sum_{\colt{y_2 \leq Y_2}{\hcf(y_2,v_2y_1)=1}}R_2(B;\ma{v},y_1,y_2)=
\vt'(\ma{v},y_1)\frac{B}{v_1^2v_2^2y_1^{2}}+
O_\ve\Big(\frac{B^{3/4+\ve}}{v_1v_2^{5/4}y_1^{3/2}}\Big),
$$
where
\begin{align*}
\vt'(\ma{v},y_1)&=-\frac{3}{\pi^2}\eta(v_2y_1^2) \prod_{p\mid  v_1 v_2 }
\Big(1-\frac{\chi(p)}{p}\Big)
\prod_{p\mid  v_1 v_2y_1 } \Big(1+\frac{1}{p}\Big)^{-1}\\
&\quad\times
\sum_{k_0 \mid v_1v_2y_1}
\mu(k_0)
\int_0^{1}\int_0^{1}\Big\{\frac{u^{1/4}v_1v_2y_1}{k_0v^{1/2}}\Big\}
\frac{\d u \d v}{\sqrt{1-u}}.
\end{align*}
\end{lem}

\begin{proof}
Now it is plain from \eqref{R_2} that
$$
\sum_{\colt{y_2 \leq Y_2}{\hcf(y_2,v_2y_1)=1}}
\hspace{-0.6cm}
R_2(B;\ma{v},y_1,y_2)=
\vt(\ma{v},y_1,1)
\frac{B^{1/2}}{v_2^{1/2}y_1}\sum_{k_0\mid v_1v_2y_1}
\hspace{-0.2cm}
\mu(k_0)
\sum_{y_2 \leq Y_2}\varpi(y_2)h\Big(\frac{y_2}{Y_2}\Big),
$$
where
$$
\varpi(y_2)=
\left\{
\begin{array}{ll}
\vt(\ma{v},y_1,y_2)/\vt(\ma{v},y_1,1), & \mbox{if
    $\hcf(y_2,v_2y_1)=1$,}\\
0, & \mbox{otherwise,}
\end{array}
\right.
$$
and
$$
h(v)
=-\frac{1}{2}\int_0^{1}\Big\{\frac{u^{1/4}v_1v_2y_1}{k_0v^{1/2}}\Big\}
\frac{\d
u}{\sqrt{1-u}}.
$$

From \eqref{calcultheta}, it follows that
$$
\varpi(y_2)= \frac{\phi(y_2)}{y_2}
\prod_{ p \mid \hcf(v_1,y_2) }\Big(1-\frac{1+\chi(p)}{p}\Big)^{-1},
$$
if  $\hcf(y_2,v_2y_1)=1$.
In order to estimate $\sum_{y_2\leq Y_2}\varpi(y_2)h(y_2/Y_2)$, we must
first calculate
the corresponding Dirichlet series $F(s)=\sum_{n=1}^\infty \varpi(n)n^{-s}$.
Let
$$
\ve_p=\Big(1-\frac{1+\chi(p)}{p}\Big)^{-1},
$$
for any prime $p$.
Then for $\Re e (s)>1$ we see that
$$
F(s)
=\prod_{p\nmid v_1v_2y_1}\Big(1+\frac{1-1/p}{p^s-1}\Big)
\prod_{p\mid v_1, p\nmid
    v_2y_1}\Big(1+\frac{\ve_p(1-1/p)}{p^s-1}\Big)
=\zeta(s)G(s),
$$
where
$$
G(s)= \prod_{ p }\Big(1-\frac{1}{p^{s+1}}\Big) \prod_{p\mid
      v_1v_2y_1} \Big(1+\frac{1-1/p}{p^s-1}\Big)^{-1}
\prod_{p\mid v_1, p\nmid v_2y_1}\Big(1+\frac{\ve_p(1-1/p)}{p^s-1}\Big).
$$
In particular it is clear that $F(s)/\zeta(s)$ is holomorphic for $\Re
e(s)>0$. It therefore follows from a standard argument (see
\cite[Lemma 2]{1} for example) that
$$
\sum_{y_2\leq t} \varpi(y_2)=  G(1)t + O(G^+(1/2)t^{1/2}),
$$
for any $t\geq 1$,
where if $G(s)=\sum_{n=1}^\infty g(n)n^{-s}$ then
$G^+(s)=\sum_{n=1}^\infty |g(n)|n^{-s}$.  A simple
calculation reveals that $G^+(1/2)\ll 2^{\omega(v_1v_2y_1)}\ll_\ve B^{\ve}$,
and so an application of partial summation yields
\begin{align*}
\sum_{y_2 \leq Y_2}&\varpi(y_2)h(y_2/Y_2)\\
   &= -\frac{G(1)}{2}\int_0^{Y_2}\!\!
\int_0^{1}\Big\{\frac{u^{1/4}v_1v_2y_1 }{k_0(s/Y_2)^{1/2}}\Big\}
\frac{\d u \d s}{\sqrt{1-u}}
+O_\ve(B^\ve Y_2^{1/2})\\
   &=
-\frac{G(1)}{2}Y_2\int_0^{1}\int_0^{1}\Big\{\frac{u^{1/4}v_1v_2y_1}{k_0v^{1/2}}\Big\}
\frac{\d u \d v}{\sqrt{1-u}} +O_\ve(B^{\ve}Y_2^{1/2}).
\end{align*}
Since  $\vt(\ma{v},y_1,1)\ll 2^{\omega(v_1y_1)}4^{\omega(v_2)}\ll_\ve
B^\ve$, this therefore completes the proof of Lemma \ref{sumR2}.
\end{proof}

While the precise value of $\vt'(\ma{v},y_1)$ in Lemma \ref{sumR2} is
perhaps unimportant, we will need the observation that
$$
\vt'(\ma{v},y_1) \ll 8^{\omega(v_1v_2y_1)} \ll_\ve
(v_1v_2y_1)^{\ve},
$$
for any $\ve>0$.
We are now ready to sum $R_2(B;\ma{v},y_1,y_2)$ over
all $(\ma{v},y_1) \in\N^3$ satisfying
$|\mu(v_2)|=1$ and the inequality
$v_1^4v_2^3y_1^2\leq B$, that follows from \eqref{ineg-3}.
Thus Lemma \ref{sumR2} implies that for any $\ve>0$ we have
\begin{align*}
\sum_{\colt{v_1,v_2,y_1\in\N}{v_1^4v_2^3y_1^2\leq B}}
\hspace{-0.1cm}|\mu(v_2)|
R_2(B;\ma{v},y_1,y_2)
&=B \sum_{\colt{v_1,v_2,y_1\in\N}{v_1^4v_2^3y_1^2\leq B}}
\frac{|\mu(v_2)|\vt'(\ma{v},y_1) }{v_1^2v_2^2y_1^{2}} 
+ O_\ve(B^{3/4+\ve})\\
&=\beta B +O_\ve(B^{3/4+\ve}),
\end{align*}
with
\begin{equation}\lab{beta}
\beta=\sum_{v_1,v_2,y_1 \in \N}
\frac{|\mu(v_2)|\vt'(\ma{v},y_1)}{v_1^2v_2^2y_1^{2}}.
\end{equation}

Let $(\ma{v},y_1,y_2) \in \N^4$ satisfy \eqref{ineg-3}, and
recall the definition \eqref{calcultheta} of $\vt(\ma{v},y_1,y_2)$.
We define $\varphi(\ma{v},y_1,y_2)$ to be zero if \eqref{polo-7} fails to
hold and
\begin{equation}\lab{PHI}
\varphi(\ma{v},y_1,y_2)=
\frac{\phi(v_1v_2y_1)}{v_1v_2y_1}\vt(\ma{v},y_1,y_2)
\end{equation}
otherwise.  Furthermore, let
\begin{equation}\lab{E}
R(B;\ma{v},y_1,y_2)=R_1(B;\ma{v},y_1,y_2)+R_2(B;\ma{v},y_1,y_2),
\end{equation}
where $R_1(B;\ma{v},y_1,y_2)$ is given by \eqref{R_1} and
$R_2(B;\ma{v},y_1,y_2)$ is given by \eqref{R_2}.
Then we have proved the following result.

\begin{lem}\lab{Sumy3y4}
Let $\ve>0$.  Then for any $B\geq 1$ and any $(\ma{v},y_1,y_2) \in
\N^4$ satisfying~\eqref{ineg-3}, we have
$$
S(\ma{v},y_1,y_2)= c\varphi(\ma{v},y_1,y_2)
\frac{B^{3/4}}{v_2^{1/4}y_1^{1/2}y_2^{1/2}}
   + R(B;\ma{v},y_1,y_2),
$$
where $c$ is given by \eqref{TAU},
$\varphi(\ma{v},y_1,y_2)$ is given by \eqref{PHI}, and
$R(B;\ma{v},y_1,y_2)$ is given by \eqref{E}
and satisfies
$$
\sum_{\colt{v_1,v_2,y_1,y_2\in\N}{\mbox{\scriptsize{\eqref{polo-7},\eqref{ineg-3}
hold}}}} R(B;\ma{v},y_1,y_2) =\beta B +O_\ve(B^{17/20+\ve}),
$$
where $\beta$ is given by \eqref{beta}.
\end{lem}

\subsection{Summation over the remaining variables}\lab{sublast}

In this section we complete our preliminary estimate for $N_{U,H}(B)$,
for which we first recall the definition \eqref{india} of the counting function
$N(Q_1,Q_2;B)$. Our task is to sum the main term in Lemma \ref{Sumy3y4} over all
$\ma{v},y_1,y_2$ satisfying  \eqref{polo-7} and \eqref{ineg-3}.
Define the arithmetic function
\begin{equation}\lab{Delta}
\Delta(n)=
\sum_{\colt{\ma{v},y_1,y_2\in\N}{v_1^4v_2^3y_1^2y_2^2=n}}
\frac{\varphi(\ma{v},y_1,y_2)}{v_2^{1/4}y_1^{1/2}y_2^{1/2}},
\end{equation}
for any $n \in \N$.  Then it follows from Lemma \ref{Sumy3y4} that
\begin{align*}
N(Q_1,Q_2;B)
&= c B^{3/4}
\sum_{\colt{v_1,v_2,y_1,y_2\in\N}{v_1^4v_2^3y_1^2y_2^2\leq B}}
\frac{\varphi(v_1,v_2,y_1,y_2)}{v_2^{1/4}y_1^{1/2}y_2^{1/2}} +\beta B
+O_\ve(B^{17/20+\ve})
\\
&=cB^{3/4}
\sum_{n \leq B} \Delta(n)+\beta B +O_\ve(B^{17/20+\ve}),
\end{align*}
for any $\ve>0$.
On inserting this estimate into Lemma \ref{Reduc1} we therefore obtain
the following result.

\begin{lem}\lab{Sum-all'}
Let $\ve>0$.  Then for any $B \geq 1$ we have
$$
N_{U,H}(B)=
4 cB^{3/4} \sum_{n \leq B} \Delta(n) +
\Big(\frac{12}{\pi^2}+ 4\beta\Big)B +  O_\ve(B^{17/20+\ve}),
$$
where $c$ is given by \eqref{TAU},
$\D(n)$ is given by \eqref{Delta} and $\beta$ is given by \eqref{beta}.
\end{lem}

\section{The height zeta function}\lab{height}

In this section we complete the proof of Theorem \ref{main'}.
Recall the definition \eqref{Delta} of $\D(n)$, which has corresponding
Dirichlet series
$$
D(s)=\sum_{n=1}^\infty \frac{\Delta(n)}{n^s}
=\sum_{v_1,v_2,y_1,y_2 \in \N}
\frac{\varphi(v_1,v_2,y_1,y_2)}{v_1^{4s}v_2^{3s+1/4}y_1^{2s+1/2}y_2^{2s+1/2}},
$$
say, where $\varphi(v_1,v_2,y_1,y_2)$ is given by \eqref{PHI}.
For $\Re e(s)>1$ the height zeta function is given by
$$
Z_{U,H}(s)=\sum_{x \in U(\Q)}\frac{1}{H(x)^s}=s\int_1^\infty
t^{-s-1}N_{U,H}(t)\d t.
$$
Thus it follows from Lemma \ref{Sum-all'} that 
$Z_{U,H}(s)=Z_1(s)+Z_2(s)$, where
\begin{align*}
Z_1(s)&=
4cs\int_1^\infty t^{-s-1/4}\sum_{n\leq t}\Delta(n) \d t
=\frac{16c\, s}{4s-3}D(s-3/4),\\
Z_2(s)&=\frac{{12/\pi^2}+4\be}{s-1}
+G_2(s),
\end{align*}
and 
$$
G_2(s)=s\int_{1}^\infty t^{-s-1}R(t)\d t 
$$ 
for some function $R(t)$ such that $R(t)\ll_\varepsilon
t^{17/20+\varepsilon}$ for any $\varepsilon>0$.  But then it follows that 
$G_2(s)$ is holomorphic on the half-plane $\{ s\in \C :\Re e(s)\geq
17/20+\varepsilon\},$ and is easily seen to satisfy the inequality 
$$
G_2(s)\ll_\ve (1+|\Im m (s)|)^{20\max\{1-\Re e (s), 0\}/3+\varepsilon}. 
$$
on this domain, via the  Phragm\'en-Lindel\"of Theorem.

It remains to analyse the function $Z_1(s)$.
The multiple sum $D(s)$ can be written as an Eulerian product
$\prod_p D_p(s)$, say.
When $p$ is odd  a routine calculation reveals that
\begin{align*}
D_p(s+1/4)&=1+\frac{(1-1/p)(2+\chi
(p))}{p^{1 +2s}-1}\Big(1+\frac{1-1/p}{p^{1+4s}-1}\Big)\\&
+\frac{(1-1/p)(1-(1+\chi (p))/p) }{p^{1+4s}-1}
   + \frac{(1-1/p)^2(1+\chi (p))}{p^{1 +3s}(1-p^{-1-4s})(1-p^{-1 -2s})}.
\end{align*}
Similarly, in the case $p=2$ we have
$$
D_2(s+1/4)=1+ \frac{1}{ 2^{1 +2s}-1}\Big(\frac{1}{2}+\frac{1}{
4(2^{1+4s}-1)}\Big)+\frac{1}{4( 2^{1+4s}-1)}\Big(1+\frac{1}{
    2^{ 3s} }\Big)+\frac{1}{2^{2+3s}} .
$$
Recall the definition \eqref{e1} of $E_1(s+1)$.  We have therefore
shown that there exists a function $H$, which is analytic on the half
plane  $\Re e(s)>1/8$, such that
\begin{equation}\lab{DE1H}
D(s+1/4)= E_1(s+1)H(s).
\end{equation}

Let $H_p(s)$ denote the Eulerian factor of $H(s)$ at the prime $p$.
We examine the analytic properties of $H(s)$ by estimating the size of
$H_p(s)$. In doing so we restrict ourselves to the domain
${\mcal{D}}=\{ s\in \C  : \Re e(s)\geq
-1/4+\varepsilon\}$.
When $p$ is odd we plainly have
$$
D_p(s+1/4)(1-p^{-1-4s})(1-p^{-1 -2s}) =1+\frac{  1+\chi
(p) }{ p^{1 +2s} } + \frac{  1+\chi (p) }{p^{1 +3s}  }
+O\Big( \frac{1 }{p^{1+\varepsilon}}\Big).
$$
It therefore follows from \eqref{DE1H} that
$$
\frac{H_p(s)}{(1-p^{-1 -3s})(1-\chi(p)p^{-1 -3s})}
=1 +\frac{
1+\chi (p) }{p^{1 +3s}}
-\frac{ 2( 1+\chi (p)) }{p^{2+5s}  } +O\Big( \frac{1 }{p^{1+\varepsilon}}\Big),
$$
whence
\begin{align*}
H_p(s) &= 1 -\frac{  2+\chi (p) }{p^{2+6s}
}-\frac{ 2( 1+\chi (p)) }{p^{2+5s}  }+\frac{  1+\chi (p) }{p^{3+9s}}
+O\Big( \frac{1}{p^{1+\varepsilon}}\Big)\\
&= E_{2,p}(s+1)\Big(1+ O\Big( \frac{1}{p^{1+\varepsilon}}\Big)\Big)
\end{align*}
on $\mcal{D}$.  Here $E_{2,p}(s+1)$ is the Eulerian factor of
$E_2(s+1)$, as given by \eqref{e2}.
Define the function
\begin{equation}\lab{calculG1}
G_1(s)=\frac{16c s}{4s-3}\frac{D(s-3/4)}{E_1(s)E_2(s)} =
\frac{16c s}{4s-3}\frac{H(s-1)}{ E_2(s)}.
\end{equation}
Then it follows from our work that $G_1(s)$ is analytic and bounded on
the half-plane $\{s\in\C: \Re e(s)\geq 3/4+\varepsilon\}$.
In order to complete the proof of Theorem \ref{main'} we simply note
from \eqref{DE1H} that
\begin{equation}\lab{H0}
H(0)= \frac{5}{2^5}\prod_{ p>2
}\Big(1-\frac{1}{p}\Big)^{4 }\Big(1-\frac{\chi(p)}{p }\Big)^2
\Big(1+\frac{4+2\chi(p)}{p}+
\frac{1}{p^2}\Big)=\tau,
\end{equation}
by \eqref{deftauf}, whence
$$
G_1(1)=16c\frac{H(0)}{E_2(1)}\neq 0.
$$

\section{Deduction of Theorem \ref{main}}\lab{deduc}

This section is virtually identical to the corresponding argument in
\cite[\S 7]{1}, and so we will be brief.
Let $\ve>0$ and let $T\in[1,B]$.
Then it follows from an application of Theorem \ref{main'}, Lemma \ref{Sum-all'}
and Perron's formula that 
\begin{align}\begin{split}\lab{pint}
N_{U,H}(B)-
\Big(\frac{12}{\pi^2}+4\beta\Big) B 
=& \frac{1}{2\pi i}
\int_{1+\ve-iT}^{1+\ve+iT} E_1(s)E_2(s)G_1(s)\frac{B^s}{s}\d s 
\\&+
O_\ve\Big(\frac{B^{1+\ve}}{T}+B^{17/20+\ve}\Big).
\end{split}\end{align}
We apply Cauchy's residue theorem to the rectangular contour $\mcal{C}$
joining the points ${\kappa-iT}$, ${\kappa+iT}$,
${1+\ve+iT}$ and ${1+\ve-iT}$, for any $\kappa \in [29/32,1)$.
On expanding the product of zeta functions about
$s=1$ it follows from \eqref{e1} that 
$$
E_1(s)=\frac{1}{48(s-1)^{4}}+O\Big(\frac{1}{(s-1 )^{3}}\Big),
$$
from which it is easy to deduce from \eqref{calculG1} and
\eqref{H0} that
$$
\mathrm{Res}_{s=1}\left\{ E_1(s)E_2(s)G_1(s)\frac{B^s}{s}\right\}=
\frac{c}{3!\times 4}L(1,\chi)^2\tau B Q(\log B),
$$
for some monic polynomial $Q$ of degree $3$.

Define the difference
$$
E(B)=N_{U,H}(B)-\frac{c}{3!\times 4}L(1,\chi)^2\tau B Q(\log B)- 
\Big(\frac{12}{\pi^2}+4\beta\Big) B.
$$
Then, in view of \eqref{pint} and the fact that the function
$E_2(s)G_1(s)$ is holomorphic and bounded for $\Re 
e (s)>5/6$, we deduce that
\begin{align}\begin{split}
E(B)&\ll_\ve
\frac{B^{1+\ve}}{T}+B^{17/20+\ve}+
\Big(\int_{\kappa-iT}^{\kappa+iT}+\int_{\kappa-iT}^{1+\ve-iT}+
\int_{1+\ve+iT}^{\kappa+iT}\Big)
\Big|E_1(s)\frac{B^s}{s}\Big|\d s,
\end{split}\lab{t3}
\end{align}
for any $\kappa \in [29/32,1)$ and any $T \in [1,B]$.
We begin by estimating the contribution from the horizontal contours.
Recall the well-known convexity bounds
$$
\zeta(\sigma+i t)\ll_\ve |t|^{(1-\sigma)/3+\varepsilon},
\quad
L(\sigma+i t,\chi)\ll_\ve
|t|^{1-\sigma +\varepsilon},
$$
that are valid for any $\sigma\in[1/2,1]$ and $|t|\geq 1$.  Then it 
follows that
$$
E_1(\sigma+it)\ll_\varepsilon
|t|^{26(1-\sigma)/3+\varepsilon}
$$
for any $\sigma\in[29/32,1]$  and $|t|\geq 1$.
This estimate allows us to deduce that
\begin{equation}
\Big(\int_{\kappa-iT}^{1+\ve-iT}+
\int_{1+\ve+iT}^{\kappa+iT}\Big)
\Big|E_1(s)\frac{B^s}{s}\Big|\d s
\ll_\ve \frac{B^{1+\ve}T^\ve}{T} + B^{\kappa }T^{\ve}.
\lab{t2}
\end{equation}

We now turn to the size of the integral
\begin{equation}\lab{t1}
\int_{\kappa-iT}^{\kappa+iT}\Big|E_1(s)\frac{B^s}{s}\Big|\d s \ll
B^\kappa\int_{-T}^{T} \frac{|E_1(\kappa+it)|}{1+|t|}  \d t =
B^\kappa I(T),
\end{equation}
say.  For given $0<U \ll T$, we will estimate the contribution
to $I(T)$ from each integral
$$
\int_{U}^{2U} \frac{|E_1(\kappa+it)|}{1+|t|}\d t \ll
\frac{1}{U}\int_{U}^{2U}|E_1(\kappa+it)|\d t = \frac{J(U)}{U},
$$
say. For any $\ve>0$ it follows from a result due to Heath-Brown
\cite{hb} that
\begin{equation}\lab{meanzeta}
\int_{U}^{2U}|\zeta(\sigma+i t)|^{8} \d t \ll_\ve U^{1+\ve},
\end{equation}
for any $\sigma \in [5/8,1]$, and any $U \geq 1$.
Moreover we will need the upper bound
\begin{equation}\lab{meanLchi}
\int_{U}^{2U}|L(\sigma+i t,\chi)|^{4} \d t \ll_\ve U^{1+\ve},
\end{equation}
that is valid for any $\sigma\in[1/2,1]$ and any $U \geq 1$.
This follows on combining standard convexity estimates
\cite[\S 7.8]{t} with work of Montgomery \cite[Theorem 10.1]{Montgomery}.
Returning to our estimate for $J(U)$, for fixed $0<U \ll T$ and any
$\kappa \in [29/32,1)$, we apply H\"older's inequality to deduce that
$$
J(U) \leq J_1^{1/8}J_2^{1/8}J_3^{1/4}J_4^{1/4}J_5^{1/4},
$$
where
\begin{align*}
J_1&=\int_{U}^{2U}|\zeta(4\kappa-3  +4it)|^{8}\d t, \qquad
J_2=\int_{U}^{2U}|\zeta(3\kappa-2+3it)|^{8}\d t,\\
J_3&=\int_{U}^{2U}|\zeta(2\kappa-1+2it)|^{8}\d t, \qquad
J_4=\int_{U}^{2U}|L(3\kappa-2+3it,\chi)|^{4}\d t,\\
J_5&=\int_{U}^{2U}|L(2\kappa-1+2it,\chi)|^{4}\d t.
\end{align*}
It therefore follows from \eqref{meanzeta} and \eqref{meanLchi} that
$
J(U)\ll_\ve U^{1+\ve},
$
on re-defining the choice of $\ve$.
Summing over dyadic
intervals for $0<U \ll T$ we obtain
$$
\int_{0}^{T}\frac{|E_1(\kappa+it)|}{1+|t|}\d t \ll_\ve T^\ve.
$$
Since we obtain the same estimate for the integral over the interval
$[-T,0]$, we deduce that $I(T)\ll_\ve T^\ve$.  We may
insert this estimate into \eqref{t1}, and then combine it with
\eqref{t2} in \eqref{t3}, in order to conclude that
$$
E(B)\ll_\ve \frac{B^{1+\ve}}{T} +B^{\kappa+\ve},
$$
for any $T\in [1,B]$.
We therefore complete the proof of Theorem \ref{main} by taking
$T=B$.


\begin{thebibliography}{99}


\bibitem{b-t'} V.V. Batyrev and Y. Tschinkel, Manin's conjecture for
   toric varieties.
\emph{J. Alg. Geom.} {\bf 7} (1998),  15--53.


\bibitem{1} R. de la Bret\`eche and T.D. Browning,
On Manin's conjecture for singular del Pezzo surfaces of degree four,
{I}. {\em Michigan Math. J.} {\bf 55} (2007), 51--80.

\bibitem{gauss} T.D. Browning,
An overview of Manin's conjecture for del Pezzo surfaces. 
{\em Proceedings of the Gauss--Dirichlet conference (G\"ottingen)}, to appear.


\bibitem{ct}
A. Chambert-Loir and Y. Tschinkel,
On the distribution of points of bounded height 
on equivariant compactifications of vector groups.
{\em  Invent. Math.} {\bf 148} (2002), 421--452.


\bibitem{c-t}
D.F. Coray and M.A. Tsfasman,
Arithmetic on singular Del Pezzo surfaces.
{\em Proc. London Math. Soc.} {\bf 57} (1988), 25--87.

\bibitem{d-t}
U. Derenthal and Y. Tschinkel,
Universal torsors over del Pezzo surfaces and rational points.
{\em Equidistribution in Number Theory: An Introduction (Montreal, 2005)}, 169--196,
NATO Science Series II: Math, Physics and Chemistry {\bf 237},
Springer, 2006.

\bibitem{f-m-t}
J. Franke, Y.I. Manin and Y. Tschinkel,
Rational points of bounded height on {F}ano varieties.
{\em Invent. Math.} {\bf 95} (1989), 421--435.

\bibitem{hb}
D.R. Heath-Brown,
Mean values of the zeta-function and divisor problems.
{\em Recent progress in analytic number theory}, Vol. I, 115--119,
Academic Press, 1981.

\bibitem{weyl}
D.R. Heath-Brown,
Weyl's inequality, Hua's inequality, and Waring's problem.  {\em
J. London Math. Soc.} {\bf 38}  (1988),  no. 2, 216--230.


\bibitem{ht}
B. Hassett and Y. Tschinkel,
Universal torsors and Cox rings.
{\em Arithmetic of higher-dimensional algebraic varieties (Palo Alto,
2002)}, 149--173,
Progr. Math. {\bf 226}, Birkh\"auser, 2004.


\bibitem{h-p} W.V.D. Hodge and D. Pedoe, \emph{Methods of algebraic
    geometry}. Vol. 2, Cambridge University Press, 1952. 


\bibitem{lipman}
J. Lipman, 
Rational singularities, with applications to algebraic surfaces and unique factorization.
{\em Inst. Hautes \'Etudes Sci. Publ. Math.} {\bf 36} (1969), 195--279.



\bibitem{Montgomery} H.L Montgomery, {\em Topics in multiplicative
number theory}. Springer Lecture Notes 227, Springer, 1971.

\bibitem{p}
E. Peyre,
Hauteurs et mesures 
de Tamagawa sur les vari\'et\'es de Fano.
{\em Duke Math. J.} {\bf 79} (1995), 101--218.

\bibitem{swd}
P. Swinnerton-Dyer,
Counting points on cubic surfaces, {II}.
\emph{Geometric methods in algebra and number theory},  
303--310,
Progr. Math. {\bf 235}, Birkh\"auser, 2005.

\bibitem{t}
E.C. Titchmarsh, {\em The theory of the Riemann zeta-function}. 2nd
ed., edited by D.R. Heath-Brown,
Oxford University Press, 1986.

\bibitem{vaaler}
J.D. Vaaler,
Some extremal functions in Fourier analysis.
{\em Bull. Amer. Math. Soc.}  {\bf 12}  (1985),  no. 2, 183--216.


\end{thebibliography}
\end{document}